\newcommand{\Lie}{{\cal L}}

\newcommand{\su}{{\frak su}}
\newcommand{\pbar}{{\overline{\p}}}

\newcommand{\diag}{{\rm diag}}

\newcommand{\tr}{{\rm tr}}

\newcommand{\F}{{\cal F}}

\newcommand{\p}{\partial}

\newcounter{note}

\newcounter{notelist}

\newtheorem{definition}{Definition}
\newtheorem{theorem}{Theorem}
\newtheorem{lemma}[theorem]{Lemma}
\newtheorem{proposition}[theorem]{Proposition}

\newtheorem{corollary}[theorem]{Corollary}

\def\d{\hbox{\rm d}}

\def\hide#1{}

\documentclass{article}
\usepackage{amsfonts}
\usepackage{amssymb}
\usepackage{amsbsy}

\begin{document}

\begin{flushright}
\end{flushright}
\vspace{0.2cm}
\begin{center}
\begin{Large}
\textbf{Surfaces associated with theta function solutions of the periodic 2D-Toda lattice}
\end{Large}\\
\vspace{1.0cm} A. M. Grundland$^{\sharp\dagger}$\footnote{
This work was supported in part by research grants from NSERC of
Canada. M. Y. Mo acknowledges the postdoctoral fellowship awarded
by the laboratory of mathematical physics of the CRM, universit\'e
de Montreal.}, M. Y.
Mo$^{\sharp}$
\bigskip\\
$^\sharp$ Centre de recherche math\'ematiques, Universit\'e\ de
Montr\'eal, C. P. 6128, Succursale centre-ville, Montr\'eal
(Quebec) H3C 3J7, Canada
\bigskip\\
$^\dagger$ Universit\'e du Qu\'ebec, Trois-Rivi\'eres CP500 (QC)
G9A 5H7, Canada
\bigskip\\
Keywords: Surface immersions, Toda lattice
\bigskip\\
\end{center}

\begin{abstract}
The objective of this paper is to present some geometric aspects
of surfaces associated with theta function solutions of the
periodic 2D-Toda lattice. For this purpose we identify the
$(N^2-1)$-dimensional Euclidean space with the ${\frak su}(N)$
algebra which allows us to construct the generalized Weierstrass
formula for immersion for such surfaces. The elements
characterizing surface like its moving frame, the Gauss-Weingarten
and the Gauss-Codazzi-Ricci equations, the Gaussian curvature, the
mean curvature vector and the Wilmore functional of a surface are
expressed explicitly in terms of any theta function solution of
the Toda lattice model. We have shown that these surfaces are all mapped
into subsets of a hypersphere in $\mathbb{R}^{N^2-1}$.
A detailed implementations of the obtained
results are presented for surfaces immersed in the ${\frak su}(2)$
algebra and we show that different Toda lattice data correspond to
different subsets of a sphere in $\mathbb{R}^3$.
\end{abstract}

\section{Introduction}\setcounter{equation}{0}\setcounter{theorem}{0} \setcounter{definition}{0}
\setcounter{example}{0}
Modern surface theory is a subject that has generated a great deal
of interests and activities in several branches of mathematics as
well as in various areas of physical science. In particular,
surfaces associated to $\mathbb{CP}^N$ sigma models provide us
with a rich class of geometric objects (See, e.g. Helein
\cite{He01}, \cite{He02}). Until very recently, the immersion of
2-dimensional surfaces obtained through the Weiestrass formula for
these models has been known only in low dimensional Euclidean
spaces. The expressions describing locally minimal surfaces (i.e.
has zero mean curvature) immersed in 3-dimensional Euclidean space
were first derived one and a half century ago by A. Enneper
\cite{E}, and K. Weierstrass \cite{We}. It starts by introducing a
pair of holomorphic functions $\psi_1$ and $\psi_2$ and then one
can introduce the 3-component real-valued vector
$X(\xi,\overline{\xi})$ as follows
\begin{eqnarray}\label{eq:im}
X(\xi,\overline{\xi})={\bf
Re}\left(i\int_{\gamma}(\psi_1^2+\psi_2^2)d\xi^1,\int_{\gamma}(\psi_1^2-\psi_2^2)d\xi^1,
-2\int_{\gamma}\psi_1\psi_2d\xi^1\right)
\end{eqnarray}
where $\gamma$ is an arbitrary curve in $\mathbb{C}$. The
parametric lines ${\bf Re}\xi={\textrm{const.}}$ and ${\bf
Im}\xi={\textrm{const.}}$ on these type of surfaces are the
minimal lines. This expression proved to be the most powerful tool
for construction and investigation of this type of surfaces in
$\mathbb{R}^3$ (for a review of the subject, see e.g. Kenmotsu
\cite{Ken}).

An extension of the Weierstrass representation (\ref{eq:im}) for
surfaces immersed into multi-dimensional Riemannian spaces in
connection with linear systems (e.g. multi-dimensional Dirac
equations) has been developed over the last decades by many
authors (e.g. \cite{BR}, \cite{DPF}, \cite{F}, \cite{KT}). The main goal is to
provide a self-contained, comprehensive approach to the
Weierstrass formulae for immersions. For this purpose, it is
convenient to formulate equations defining the immersion directly
in the matrix form which take their values in a finite dimensional
Lie algebra. The main advantages of this procedure appear when
group analysis of the immersion makes it possible to construct
regular algorithms for finding certain classes of surfaces without
referring to any additional considerations. Instead, we proceed
directly from the given model of equations. The task of finding an
increasing number of surfaces is facilitated by the group
properties of the model considered. The rich character of this
formulation makes the immersion formula a rather interesting object to
investigate.

This paper is concerned with some geometric aspects of smooth
2-dimensional surfaces arising from the study of theta function
solutions of the 2D Toda lattice model. To achieve this we
identify $\mathbb{R}^{N^2-1}$ with the Lie algebra ${\frak su}(N)$
which allows us to construct the generalized Weierstrass formula
for immersion of surfaces. This formula is a consequence of
conservation laws of the model. We show that there exist $N^2-1$
real-valued functions which can be treated as coordinates of the
surface immersed in the ${\frak su}(N)$ algebra. We reformulate
explicitly the structural equations for the immersion in Cartan's
language of moving frames on a surface. These are expressed in
terms of any theta function solution of the Toda model. In fact,
this approach allows us to find explicitly the form of the
Gauss-Weingarten and the Gauss-Codazzi-Ricci equations in terms of
this type of solution. The first and second fundamental forms, the
Gauss curvature and the mean curvature vector are the central
concepts described here. We present in details an implementation
of these results for surfaces immersed in the ${\frak su}(2)$
algebra. In the case of ${\frak su}(2)$ we show that different
Toda lattice data correspond to different subsets of a sphere in
$\mathbb{R}^3$.

From the point of view of the 2D-Toda lattice model, this paper
addresses the following two questions.

1. Starting from the theta function of the 2D-Toda lattice model one
can derive the conservation laws and immersion of a 2D-surface in
multi-dimensional space $\mathbb{R}^{N^2-1}={\frak su}(N)$. What
are the properties of this surface?

2. Given that there is a connection between 2D-Toda lattice models
and surfaces immersed in ${\frak su}(N)$, do these surfaces tell
us anything new? Is this connection useful?

Finally, let us note that the outlined approach to the study of surfaces
leads to the following potential applications of the theory of
surfaces associated to the doubly periodic solutions of the
2D-Toda lattice.

In physics the doubly periodic solutions have found applications
in such varied areas as field theory, quantum field theory and
string theory \cite{A}, \cite{IP}, \cite{Po}, statistical physics
\cite{DG}, \cite{NPW}, phase transitions (e.g. growth of crystals,
deformations of membranes, dynamics of vortex sheets, surfaces
waves etc. \cite{CJZ}), fluid dynamics, e.g. the motion of
boundaries between regions of different densities and viscosities
\cite{JR}.

In biochemistry and biology, surfaces have been shown to play an
essential role in several applications to nonlinear phenomena in
the study of biological membranes and vesicles, for example long
protein molecules \cite{LS}, \cite{Da}, \cite{OLX}, the
Canham-Helfrich membrane models \cite{La}. These macroscopic
models can be derived from microscopic ones and allow us to
explain basic features and equilibrium shapes both for biological
membranes and for liquid interfaces \cite{DM}, \cite{S}. In
chemistry there are applications e.g. to energy and momentum
transport along a polymer molecule \cite{Da}.

The results in this paper give a systematic way of constructing
surfaces from doubly periodic solutions which could be applied to
the above areas of research.

In mathematics, the analysis described in this paper could be
extended to a systematic description of surfaces in the Lie algebras/
Lie groups framework and isomonodromic deformations in connection
with surfaces (e.g. \cite{B}, \cite{BE} and the references
therein).

The paper is organized as follows. In section \ref{se:DS}, we give
a brief introduction to the standard techniques in integrable
systems and in particular the 2D-Toda lattice. In section
\ref{se:theta}, we explained how to construct solutions of the
2D-Toda lattice from theta functions of Riemann surfaces, while
section \ref{se:immersion} defines immersions from these solutions
of the 2D-Toda lattice and introduced some basic results in
surface immersion theory. In section \ref{se:suN} we expressed
these immersions in terms of theta functions and derived the
moving frame, the first and second fundamental forms and the mean
curvature vector. Finally, section \ref{se:su2} sees the
application of these results to the ${\frak su}(2)$ case in which
the surfaces obtained have constant mean and Gaussian curvatures.

\section{Integrable systems structure of the 2D-Toda
lattice}\label{se:DS}\setcounter{equation}{0}\setcounter{theorem}{0} \setcounter{definition}{0}
\setcounter{example}{0}

The periodic Toda lattice is equivalent to the following
compatibility equations
\begin{eqnarray}\label{eq:toda}
\overline{\p}A&-&\p B=[A,B], \quad B=-A^{\dag}, \\
A&=&\pmatrix{\p u_0&0 &0&\ldots &0&U_{0,N} \cr
             U_{1,0}&\p u_1&0&\ldots &0&0\cr
             0&U_{2,1}&\p u_2&\ldots &0&0&\cr
             \vdots &\vdots &\vdots &\ddots &\vdots &\vdots\cr
             0&0&0&\ldots &\p u_{N-1}&0\cr
             0&0&0&\ldots &U_{N,N-1}&\p u_N\cr}
\end{eqnarray}
where $\p=\p_{\xi}$, $\overline{\p}=\p_{\overline{\xi}}$, $u_i$
are real-valued functions $u_i:\mathbb{C}\rightarrow\mathbb{R}$
such that $\sum_{i=1}^N u_i=0$ and $U_{i,j}=e^{u_i-u_j}$.
Therefore $A$ is in the ${\frak sl}(n,\mathbb{C})$ algebra.

To understand the integrable systems structure, first let us write
the zero curvature condition as one that involves a complex
parameter $\lambda$
\begin{eqnarray*}
\overline{\p}A_{\lambda}&-&\p B_{\lambda}=[A_{\lambda},B_{\lambda}], \quad B_{\lambda}=-\rho(A_{\lambda}), \quad \lambda\in\mathbb{C} \\
A_{\lambda}&=&\pmatrix{\p u_0&0 &0&\ldots &0&\lambda U_{0,N} \cr
             U_{1,0}&\p u_1&0&\ldots &0&0\cr
             0&U_{2,1}&\p u_2&\ldots &0&0&\cr
             \vdots &\vdots &\vdots &\ddots &\vdots &\vdots\cr
             0&0&0&\ldots &\p u_{N-1}&0\cr
             0&0&0&\ldots &U_{N,N-1}&\p u_N\cr}
\end{eqnarray*}
where $\rho$ is the involution on the series of matrix defined by
\begin{eqnarray}\label{eq:involalg}
\rho\left(\sum_{-\infty}^{\infty}X_i\lambda^i\right)=\sum_{-\infty}^{\infty}X_i^{\dag}\lambda^{-i},\quad
X_i\in{\frak sl}(n,\mathbb{C})
\end{eqnarray}
We can represent $A_{\lambda}$ as follows
\begin{eqnarray*}
A_{\lambda}&=&A_0+A_1\Lambda, \quad B=-A_0^{\dag}-A_1^{\dag}\Lambda^{-1}\\
\Lambda&=&\pmatrix{0&0 &0&\ldots &0& \lambda \cr
             1&0&0&\ldots &0&0\cr
             0&1&0&\ldots &0&0&\cr
             \vdots &\vdots &\vdots &\ddots &\vdots &\vdots\cr
             0&0&\ldots &1&0&0\cr
             0&0&0&\ldots &1&0\cr}, \quad \Lambda^{-1}=\pmatrix{0&1 &0&\ldots &0& 0\cr
             0&0&1&\ldots &0&0\cr
             0&0&0&\ldots &0&0&\cr
             \vdots &\vdots &\vdots &\ddots &\vdots &\vdots\cr
             0&0&\ldots &0&0&1\cr
             \lambda^{-1}&0&0&\ldots &0&0\cr}
\end{eqnarray*}
where $A_0$, $A_1$ are diagonal matrices.

One can now look for formal solutions $F$ of the linear equations
\begin{eqnarray}\label{eq:linear}
\p F=-A_{\lambda}F, \quad \overline{\p} F=-B_{\lambda}F
\end{eqnarray}
near $\lambda=0$ and $\lambda=\infty$ respectively. Since
$U_{i,j}=e^{u_i-u_j}$, the determinant of $A_1$ equals 1. It is a
standard result due to Drinfeld and Sokolov \cite{DS} that near
$\lambda=\infty$, there exists a formal solution $F_{\infty}$ of
the form
\begin{eqnarray}\label{eq:dressing}
F_{\infty}=T_{\infty}\exp(\xi\Lambda)
\end{eqnarray}
where $T_{\infty}$ is a formal power series in $\Lambda^{-1}$ of
the form
\begin{eqnarray}\label{eq:series}
T_{\infty}=D_0+D_1\Lambda^{-1}+\ldots
\end{eqnarray}
with diagonal matrices $D_i$ such that $D_0$ is invertible.
Similarly, near $\lambda=0$, there exists a formal solution $F_0$
of the form
\begin{eqnarray*}
F_{0}=T_{0}\exp(\overline{\xi}\Lambda^{-1})
\end{eqnarray*}
with similar properties. The matrices $T_0$ and $T_{\infty}$ are
called the dressing matrices.

Let $L_i$ and $\overline{L}_i$ denote the following matrices,
\begin{eqnarray}\label{eq:Li}
L_i=(T_{\infty}\Lambda^iT_{\infty}^{-1})_+, \quad
\overline{L}_i=(T_0\Lambda^{-i}T_0^{-1})_-, \quad
i=1,\ldots,\infty
\end{eqnarray}
where $+$, $-$ denotes the polynomial parts in $\Lambda$ and
$\Lambda^{-1}$ respectively. From now on, we will suppress the
$\lambda$ dependence of $L_i$.

Usual argument in integrable systems theory shows that one can
then define infinitely many commuting flows according to (See e.g.
\cite{DS})
\begin{eqnarray}\label{eq:compa}
\p_{t_i}A_{\lambda}-\p L_i=[A_{\lambda},L_i], \quad
\p_{\overline{t}_i}A_{\lambda}-\p \overline{L}_i=[A_{\lambda},\overline{L}_i] \\
\p_{t_i}B_{\lambda}-\overline{\p}L_i=[B_{\lambda},L_i], \quad
\p_{\overline{t}_i}B_{\lambda}-\overline{\p}
\overline{L}_i=[B_{\lambda},\overline{L}_i]\nonumber
\end{eqnarray}
To see how this works, consider, for example, the case of $L_i$.
First make a gauge transformation
\begin{eqnarray*}
\p+\hat{A}_{\lambda}&=&\mathbb{G}^{-1}(\p+A_{\lambda})\mathbb{G}\\
\mathbb{G}&=&\diag(e^{u_0},e^{u_1},\ldots,e^{u_{N}})
\end{eqnarray*}
so that $\hat{A}_{\lambda}$ has $\Lambda$ as its leading term.
That is
\begin{eqnarray*}
\hat{A}_{\lambda}=\Lambda+\hat{A}_0
\end{eqnarray*}
we have
\begin{eqnarray}\label{eq:b}
[\p +\hat{A}_{\lambda}, L_i]=[\p+\hat{A}_{\lambda},
(\mathbb{G}^{-1}T_{\infty}\Lambda^iT_{\infty}^{-1}\mathbb{G})-(\mathbb{G}^{-1}T_{\infty}\Lambda^iT_{\infty}^{-1}\mathbb{G})_{<0}]
\end{eqnarray}
where
$(\mathbb{G}^{-1}T_{\infty}\Lambda^iT_{\infty}^{-1}\mathbb{G})_{<0}$
is a series in negative powers of $\Lambda$ only. That is, it has
no constant term in $\Lambda$.

The series $T_{\infty}$ is the gauge that transforms
$\p+A_{\lambda}$ to $\p+\Lambda$, that is,
\begin{eqnarray*}
T_{\infty}(\p+\Lambda)T_{\infty}^{-1}=\p+A_{\lambda}
\end{eqnarray*}
we have
\begin{eqnarray*}
[\p +\hat{A}_{\lambda},
\mathbb{G}^{-1}T_{\infty}\Lambda^iT_{\infty}^{-1}\mathbb{G}]=\mathbb{G}^{-1}T_{\infty}[\p+\Lambda,
\Lambda^i]T_{\infty}^{-1}\mathbb{G}=0,
\end{eqnarray*}
therefore we have
\begin{eqnarray}\label{eq:a}
[\p +\hat{A}_{\lambda}, L_i]=[\p+\hat{A}_{\lambda},-(\mathbb{G}^{-1}T_{\infty}\Lambda^iT_{\infty}^{-1}\mathbb{G})_{<0}]
\end{eqnarray}
The left hand side is a positive series in $\Lambda$ while the
right hand side is a non-positive series in $\Lambda$, we see that
it must be a constant in $\Lambda$. According to the form of the
series expansion (\ref{eq:series}), both sides of (\ref{eq:a})
must be a diagonal matrix. We can therefore set this diagonal
matrix to be $\p_{t_i}\hat{A}_{\lambda}$. By applying the gauge
transformation $\mathbb{G}^{-1}$ to these equations, we obtain the
equations (\ref{eq:compa}).

The argument for $L_i$ and $B$ is similar, but we use the
following instead of (\ref{eq:b})
\begin{eqnarray*}
\overline{\p} F_{\infty}&=&-B_{\lambda}F_{\infty}\\
\overline{\p} T_{\infty}\exp(\xi\Lambda)&=&-B_{\lambda}T_{\infty}\exp(\xi\Lambda)\\
\overline{\p} T_{\infty}&=&-B_{\lambda}T_{\infty}
\end{eqnarray*}
All these flows commute with each other and each flow defines
an invariant of the 2D-Toda lattice.

\section{Theta function solutions of the Toda
lattice}\label{se:theta}\setcounter{equation}{0}\setcounter{theorem}{0} \setcounter{definition}{0}
\setcounter{example}{0}

In this section we present some basic facts concerning
multi-dimensional theta functions on Riemann surfaces that will provide
the tool to construct explicit solutions of the Toda lattice.
Later we will see how to construct surfaces from these solutions.

\subsection{Basic results of theta function}

Let us first remind ourselves some useful facts of theta function
that will allow us to construct explicitly meromorphic functions
and functions with essential singularities of exponential type.
For a review of the subject see e.g. \cite{Mu}.

First choose a canonical basis of cycles $\{a_i,b_i\}$ on a
Riemann surface $\Sigma$, and let $\omega_i$ be 1-forms dual to
this basis, that is
\begin{eqnarray*}
\int_{a_i}\omega_j=\delta_{ij}, \quad \int_{b_i}\omega_j=\tau_{ij}
\end{eqnarray*}
We can define a lattice $L(M)$ in $\mathbb{C}^g$ by using the
columns of the $g\times 2g$ matrix $(I_d, \Pi)$, where $I_d$ is the
$g\times g$ identity matrix and $(\Pi)_{ij}=\tau_{ij}$. The torus
$\mathbb{C}^g\backslash L(M)$ is called the Jacobian $Jac(\Sigma)$
of the Riemann surface.

The theta function associated to the Riemann surface is a function
$\Theta:\mathbb{C}^g\rightarrow \mathbb{C}$ defined by
\begin{eqnarray*}
\Theta (\overrightarrow{z})= \sum_{\overrightarrow{n}\in
\mathbb{Z}^g} {\rm e}^{i\pi \overrightarrow{n}\cdot \Pi
\overrightarrow{n} - 2i\pi \overrightarrow{z}\cdot
\overrightarrow{n}}
\end{eqnarray*}
The theta function has the following periodicity.
\begin{proposition}\label{pro:per}
Let $e^k$, $\tau^k$ be the columns of the matrices $I_d$ and $\Pi$
respectively, then
\begin{eqnarray}\label{eq:period}
\Theta (\overrightarrow{z}+e^k)&=& \Theta(\overrightarrow{z}) \nonumber\\
\Theta (\overrightarrow{z}+\tau^k)&=& \exp 2\pi
i\left(-z_k-{{\tau_{kk}}\over{2}}\right)\Theta(\overrightarrow{z})
\end{eqnarray}
\end{proposition}
Let $U:\Sigma\rightarrow Jac(\Sigma)$ be the Abel map, than the
theta function composite with the Abel map has $g$ zeros of on
$\Sigma$. In fact, let $D=\sum_{i=1}^g\gamma_i$ be a divisor of
degree $g$ on $\Sigma$, then the following multi-valued function
\begin{eqnarray*}
\Theta (U(p)-U(D)-K)
\end{eqnarray*}
has zeros at the $g$ points $\gamma_i$, where the vector
$K=(K_1,\ldots,K_g)$ is the Riemann constant
\begin{eqnarray*}
K_j={{2\pi i+\tau_{jj}}\over 2}-{1\over{2\pi i}}\sum_{l\neq
j}\int_{a_l}(\omega_l(P)\int_{P_0}^P\omega_j)
\end{eqnarray*}

\subsection{Baker functions}
We will now define the Baker functions, which would serve as the
fundamental solutions of the linear system (\ref{eq:linear}).
\begin{definition}\label{de:baker}
Let $\Sigma$ be a Riemann surface of genus $g$, and
$Q_1,\ldots,Q_l$ are $l$ points on $\Sigma$. Let $k_i^{-1}$ be
local coordinates in neighborhoods of these points
($k_i(Q_i)=\infty$), and $D$ be a divisor on
$\Sigma/(Q_1\cup\ldots\cup Q_l)$. Then a Baker l point function
$f(P)$ is a meromorphic function on $\Sigma/(Q_1\cup\ldots\cup
Q_l)$ that has pole divisor $D$ and that near $Q_i$,
$f\exp(-q_i(k_i))$ is holomorphic for some polynomial $q_i$ in
$k_i$.
\end{definition}
A useful theorem that will be used throughout the construction is
the Riemann-Roch theorem.
\begin{theorem}\label{thm:RR} (Riemann-Roch) Let $D$ be the pole
divisor of the Baker l point function $f(P)$ as in definition
\ref{de:baker} and let $d$ be its degree. If $D$ is not a special
divisor, then the Baker l point function that has pole divisor $D$
forms a linear space that has dimension $max(d-g+1,0)$, where $d$
is the degree of $D$. (See, e.g. \cite{D81}).
\end{theorem}
In particular, if the divisor $D$ is of degree $g$, then the Baker
$l$ point function with pole divisor $D$ will be unique up to a
constant factor. If the divisor $D$ is of degree $g+n$, but $n$
zeros of the Baker $l$ point function are given, than the Baker
$l$ point function can also be uniquely determined up to a
constant factor. We therefore have the following
\begin{corollary}\label{cor:uni}
Let $D$ be a non-special divisor of degree $g$, and $Q_1$, $Q_2$
be two fix points on a Riemann surface $\Sigma$ of genus $g$, then
there exist unique meromorphic functions $f_n$ on $\Sigma/(Q_1\cup
Q_2)$ that have pole divisor $D$ and such that
\begin{eqnarray}\label{eq:asyfn}
f_n(P)&=&k_1^n\left(\sum_{i=0}^{\infty}h_{n,i}^{(1)}k_1^{-i}\right)\exp(k_1\xi),\quad
P\rightarrow Q_1, \nonumber\\
f_n(P)&=&k_2^{-n}\left(1+\sum_{i=1}^{\infty}h_{n,i}^{(2)}k_2^{-i}\right)\exp(k_2\overline{\xi}),\quad
P\rightarrow Q_2,
\end{eqnarray}
where $k_i^{-1}$ are local coordinates near $Q_i$, that is,
$k_i^{-1}(Q_i)=0$.
\end{corollary}
In fact, the functions $f_n$ can be expressed in terms of the
theta function on $\Sigma$.

Let $\Omega^1$ and $\Omega^2$ be the normalized meromorphic
1-forms on $\Sigma$ (that is, $\int_{a_j}\Omega^i=0$) that are
holomorphic on $\Sigma/Q_1$ and $\Sigma/Q_2$ respectively, and
that
\begin{eqnarray*}
\Omega^1&=&dk_1+O(k_1^{-1}),\quad P\rightarrow Q_1,\\
\Omega^2&=&dk_2+O(k_2^{-1}),\quad P\rightarrow Q_2
\end{eqnarray*}
where $O(k^{-1})$ are terms of order $k^{-1}$. Let us denote the
$b$-periods of these 1-forms by the vectors $\mathbb{B}_1$, $\mathbb{B}_2$.

We can now express the functions $f_n$ in terms of theta
functions.

\begin{proposition}\label{pro:fn}
The functions $f_n$ can be expressed in terms of theta functions
as follows
\begin{eqnarray}\label{eq:fn}
f_n(P)&=&C{{\Theta(U(P)-U(D)+n(U(Q_2)-U(Q_1))+\xi
\mathbb{B}_1+\overline{\xi}\mathbb{B}_2-K)}\over
{\Theta(U(P)-U(D)-K)}}\nonumber\\
&\times&\left({{\Theta(U(P)-U(Q_2)-e)}\over
{\Theta(U(P)-U(Q_1)-e)}}\right)^n
\exp\left(\xi\int^P\Omega^1+\overline{\xi}\int^P\Omega^2\right)
\end{eqnarray}
where $\Theta$ is the theta function and $e$ is a vector of the
form
\begin{eqnarray*}
e=U(x_1)+\ldots +U(x_{g-1})+K
\end{eqnarray*}
where $x_i$ are $g-1$ arbitrary points on $\Sigma$ and $C$ is a
normalization constant.
\end{proposition}
{\bf Proof}. Since the factors $\Theta(U(P)-U(Q_1)-e)$ and
$\Theta(U(P)-U(Q_2)-e)$ has zeros at $\{x_i, Q_1\}$ and
$\{x_i,Q_2\}$ respectively, the factor
$\left({{\Theta(U(P)-U(Q_2)-e)}\over
{\Theta(U(P)-U(Q_1)-e)}}\right)^n$ on the right hand side of
(\ref{eq:fn}) has an order $n$ zero at $Q_2$ and an order $n$ pole
at $Q_1$. Similarly, the factor $\Theta(U(P)-U(D)-K)$ has zeros at
the points in $D$. Therefore the right hand side of (\ref{eq:fn})
has poles at $D$ with the asymptotic form indicated as in
(\ref{eq:asyfn}). One can now verify that the right hand side of
(\ref{eq:fn}) is in fact single-valued by using the periodicity of
the theta function (\ref{eq:period}). $\Box$

To obtain solutions of the Toda lattice from this, we need the
following lemma \cite{D81}
\begin{lemma}\label{le:real}
Suppose that there is an involution $\rho:\Sigma\rightarrow\Sigma$
that permutes $Q_1$ and $Q_2$ such that
$k_1=-\overline{\rho(k_2)}$. If the divisor $D$ that defines $f_n$
in corollary \ref{cor:uni} is such that $D+\rho(D)$ is the zero
divisor of a differential $\omega$ of third kind with simple poles
at $Q_1$ and $Q_2$, then the coefficients $h_{n,0}^{(1)}$ are
real.
\end{lemma}
{\bf Proof}. By a differential of third kind, we mean that
$\omega$ has only simple poles at $Q_1$ and $Q_2$, with residues 1
and -1 respectively.

Since $D+\rho(D)$ is the zero divisor of $\omega$, the following
differential
\begin{eqnarray*}
\tilde{\omega}=f_n(P)\overline{f}_n(\rho(P))\omega
\end{eqnarray*}
is again a differential of third kind with simple poles only at
$Q_1$ and $Q_2$. Its residues are $(-1)^nh_{n,0}^{(1)}$ and
$-(-1)^n\overline{h}_{n,0}^{(1)}$ respectively. Since the sum of
residues is zero, the lemma follows. $\Box$

If we now let $f_n$ be functions on a Riemann surface defined by a
polynomial
\begin{eqnarray}\label{eq:RS}
y^{N+1}+a_0\lambda^m+\overline{a}_0\lambda^{-m}+P(y,\lambda)=0
\end{eqnarray}
such that $m$, $N+1$ are relatively prime and that
\begin{eqnarray*}
\overline{P}(\overline{y},\overline{\lambda}^{-1})=P(y,\lambda)
\end{eqnarray*}
where $P(y,\lambda)$ is a polynomial in $y$ and $\lambda$. We can
define an the involution $\rho$ on this Riemann surface by
\begin{eqnarray}\label{eq:invol}
\rho:(y,\lambda)\mapsto (\overline{y},\overline{\lambda}^{-1})
\end{eqnarray}
Since the Riemann surface is branched at
$Q_1=(y=\infty,\lambda=\infty)$ and $Q_2=(y=\infty,\lambda=0)$, we
can choose local coordinates $k_1$ and $k_2$ to be
\begin{eqnarray}\label{eq:coord}
k_1=\lambda^{1\over {N+1}}, \quad k_2=-\lambda^{-{1\over {N+1}}}
\end{eqnarray}
We can now state the main theorem of this section
\begin{theorem}
Suppose the $f_n$ in corollary \ref{cor:uni} defined by a divisor
$D$ satisfies the condition in lemma \ref{le:real} with local
coordinates and involutions defined by (\ref{eq:coord}) and
(\ref{eq:invol}) respectively. Then if
$h_{n,0}^{(1)}h_{n+1,0}^{(1)}<0$, the functions $f_n$ satisfy the
following
\begin{eqnarray}\label{eq:comp}
\p f_n&=&-U_{n,n-1}^2f_{n-1} \nonumber\\
\overline{\p}f_n&=&f_{n+1}+2\overline{\p}u_nf_n, \quad
e^{2u_n}=-h_{n,0}^{(1)}
\end{eqnarray}
\end{theorem}
{\bf Proof}. We shall now denote $Q_1$ by $Q_{\infty}$ and $Q_2$
by $Q_0$ and change all the indices accordingly.

To prove the first equation in (\ref{eq:comp}), recall that the
functions $f_n$ has the following asymptotic behavior
\begin{eqnarray*}
f_n(P)&=&k_{\infty}^{-n}\left(\sum_{i=0}^{\infty}h_{n,i}^{(\infty)}k_1^{-i}\right)\exp(k_{\infty}\xi),\quad
P\rightarrow Q_{\infty}, \\
f_n(P)&=&k_0^{n}\left(1+\sum_{i=1}^{\infty}h_{n,i}^{(0)}k_2^{-i}\right)\exp(k_0\overline{\xi}),\quad
P\rightarrow Q_0,
\end{eqnarray*}
Suppose that $h_{n,0}^{(1)}$ are negative and let
$e^{2u_n}=-h_{n,0}^{(1)}$, then
\begin{eqnarray*}
\p
f_n+U_{n,n-1}^2f_{n-1}=k_{\infty}^{-n}O(1)\exp(k_{\infty}\xi),\quad
P\rightarrow Q_{\infty}
\end{eqnarray*}
where $O(1)$ is a term that is holomorphic at $Q_{\infty}$, and
\begin{eqnarray*}
\p
f_n+U_{n,n-1}^2f_{n-1}=k_{0}^{n-1}O(1)\exp(k_{0}\overline{\xi}),\quad
P\rightarrow Q_{0}
\end{eqnarray*}
where $O(1)$ is a term that is holomorphic at $Q_{0}$. We see that
the left hand side is a function $R$ with poles only at $D$ and
same asymptotic behavior as $f_{n-1}$ except that near
$Q_{\infty}$, it has a zero of order $n$ while $f_{n-1}$ has a
zero of order $n-1$. Therefore ${R\over f_{n-1}}=0$ at
$Q_{\infty}$. Now by the Riemann-Roch theorem (theorem \ref{thm:RR}),
we see that $R=0$.

We now prove the second equation in (\ref{eq:comp}). We see that
\begin{eqnarray*}
\overline{\p}f_n-f_{n+1}-2\overline{\p}u_nf_n=k_0^{n}O(1)\exp(k_0\overline{\xi}),\quad
P\rightarrow Q_0
\end{eqnarray*}
while near $Q_{\infty}$, it has the following asymptotics
\begin{eqnarray*}
\overline{\p}f_n-f_{n+1}-2\overline{\p}u_nf_n=k_0^{n+1}O(1)\exp(k_{\infty}\xi
),\quad P\rightarrow Q_{\infty}
\end{eqnarray*}
Since $h_{n,0}^{(1)}=-e^{2u_n}$. Then by similar argument as
before, we see that the both sides must be zero. $\Box$

We can now express the solutions in terms of theta functions.
\begin{proposition}\label{pro:sol}
The functions $u_i$ that solve the Toda lattice equation
(\ref{eq:toda}) can be expressed in terms of theta function as
\begin{eqnarray}\label{eq:sol}
u_n&=&{1\over 2}\log
\left|{{\Theta(U(Q_{\infty})-U(D)+n(U(Q_0)-U(Q_{\infty}))+\xi
\mathbb{B}_1+\overline{\xi}\mathbb{B}_2-K)}\over
{\Theta(U(Q_0)-U(D)+n(U(Q_0)-U(Q_{\infty}))+\xi
\mathbb{B}_1+\overline{\xi}\mathbb{B}_2-K)}}\right| \nonumber\\
&+&c+c_nn
\end{eqnarray}
where $c$ and $c_n$ are constants.
\end{proposition}
{\bf Proof}. We first show that the set of equations
(\ref{eq:comp}) implies (\ref{eq:toda}). To see this, note that
the condition
\begin{eqnarray*}
\overline{\p}A_{\lambda}-\p B_{\lambda}=[A_{\lambda},B_{\lambda}]
\end{eqnarray*}
is equivalent to
\begin{eqnarray*}
[\p+A_{\lambda},\overline{\p}+B_{\lambda}]=0
\end{eqnarray*}
We can make a gauge transformation $\p+A_{\lambda}\mapsto
\mathbb{G}(\p+A_{\lambda})\mathbb{G}^{-1}$, $\p+B_{\lambda}\mapsto
\mathbb{G}(\p+B_{\lambda})\mathbb{G}^{-1}$ by the gauge
\begin{eqnarray*}
\mathbb{G}=\diag(e^{u_0},e^{u_1},\ldots,e^{u_N})
\end{eqnarray*}
This will change the Lax pair into the following
\begin{eqnarray}\label{eq:aprime}
\p+\hat{A}_{\lambda}&=&\p+\pmatrix{0&0 &0&\ldots &0&\lambda
U_{0,N}^2 \cr
             U_{1,0}^2&0&0&\ldots &0&0\cr
             0&U_{2,1}&0&\ldots &0&0\cr
             \vdots &\vdots &\vdots &\ddots &\vdots &\vdots\cr
             0&0&0&\ldots &0&0\cr
             0&0&0&\ldots &U_{N,N-1}^2&0\cr} \nonumber\\
\p+\hat{B}_{\lambda}&=&\p+\pmatrix{2\overline{\p}u_0&1 &0&\ldots
&0&0 \cr
             0&2\overline{\p}u_1&1&\ldots &0&0\cr
             0&0&2\overline{\p}u_2&1&\ldots &0\cr
             \vdots &\vdots &\vdots &\ddots &\vdots &\vdots\cr
             0&0&0&\ldots &2\overline{\p}u_{N-1}&1\cr
             \lambda^{-1}&0&0&\ldots &0&2\overline{\p}u_N\cr}
\end{eqnarray}
Let $f_n$ be functions on a Riemann surface $\Sigma$ (\ref{eq:RS})
defined by (\ref{eq:fn}). The Riemann surface forms a $N+1$ sheet
covering of the $\lambda$-plane. Away from the branch points, each
point $P=\lambda$ on the $\lambda$-plane corresponds to $N+1$
points $P_i=(y_i,\lambda)$ on the Riemann surface $\Sigma$. Since
the $f_n$ provides a solution to (\ref{eq:comp}), if we let $F_1$
be the following matrix
\begin{eqnarray}\label{eq:F1}
F_1=\pmatrix{f_0(P_1)&f_0(P_2)&\ldots &f_0(P_{N+1}) \cr
             f_1(P_1)&f_1(P_2)&\ldots &f_1(P_{N+1})\cr
             f_2(P_1)&f_2(P_2)&\ldots &f_2(P_{N+1})\cr
             \vdots &\vdots &\vdots &\vdots \cr
             f_{N-1}(P_1)&f_{N-1}(P_2)&\ldots &f_{N-1}(P_{N+1})\cr
             f_{N}(P_1)&f_{N}(P_2)&\ldots &f_{N}(P_{N+1})\cr}
\end{eqnarray}
then the following equations will be satisfied
\begin{eqnarray}\label{eq:hat}
(\p+\hat{A}_{\lambda})F_1=0, \quad (\p+\hat{B}_{\lambda})F_1=0
\end{eqnarray}
Most of the equations in the above follows from (\ref{eq:comp}),
The only non trivial ones are the ones that involve $\lambda$,
that is
\begin{eqnarray*}
\p f_0&=&-\lambda U_{0,N}^2f_N \\
\overline{\p} f_N&=&\lambda^{-1}f_0+\p u_Nf_N
\end{eqnarray*}
To see that this is true, we observe that $\lambda^{-1}$ is a
meromorphic function on $\Sigma$ that has an order $N$ pole at the
point $Q_0$ and an order $N$ zero at the point $Q_{\infty}$ with
no poles or zeros elsewhere. This means that
$\lambda^{-1}f_0=f_{N+1}$ and hence the above equations are true.

Since the hatted ($\hat{A}_{\lambda}$) connection (\ref{eq:hat}) and the unhatted one is related by the
gauge transformation $\mathbb{G}$, a fundamental solution $F$ to
the equations (\ref{eq:linear}) can now be formed by using the
$f_n$
\begin{eqnarray*}
F(\lambda)=\mathbb{G}^{-1}F_1
\end{eqnarray*}
It remains to express $h_{n,0}^{(1)}$ in terms of theta functions.
Since $\Theta(U(P)-U(Q_{\infty})-e)$ and $\Theta(U(P)-U(Q_0)-e)$
have zero at $Q_{\infty}$ and $Q_0$ respectively, we have
\begin{eqnarray*}
\Theta(U(P)-U(Q_{\infty})-e)&=&s_{\infty}k_{\infty}^{-1}+\ldots,\quad P\rightarrow Q_{\infty}\\
\Theta(U(P)-U(Q_{0})-e)&=&s_{0}k_{0}^{-1}+\ldots,\quad P\rightarrow Q_0\\
\end{eqnarray*}
So the normalization constant in (\ref{eq:fn}) is
\begin{eqnarray*}
C&=&\left({{\Theta(U(Q_0)-U(Q_{\infty})-e)^n}\over
{\Theta(U(Q_0)-U(D)-K)s_0^n}}\right)^{-1}\\
&\times&\Theta(U(Q_{0})-U(D)+n(U(Q_0)-U(Q_{\infty}))+\xi
\mathbb{B}_1+\overline{\xi}\mathbb{B}_2-K)
\end{eqnarray*}
therefore $h_{n,0}^{(1)}$ is
\begin{eqnarray*}
h_{n,0}^{(1)}&=&
{{\Theta(U(Q_{\infty})-U(D)+n(U(Q_0)-U(Q_{\infty}))+\xi
\mathbb{B}_1+\overline{\xi}\mathbb{B}_2-K)}\over
{\Theta(U(Q_0)-U(D)+n(U(Q_0)-U(Q_{\infty}))+\xi
\mathbb{B}_1+\overline{\xi}\mathbb{B}_2-K)}}e^ce^c_nn
\end{eqnarray*}
where $c$ and $c_n$ are constants
\begin{eqnarray*}
c&=&\ln{{\Theta(U(Q_0)-U(D)-K)}\over{\Theta(U(Q_{\infty})-U(D)-K)}}\\
 c_n&=&\ln{{\Theta
(U(Q_{\infty})-U(Q_0)-e)s_0}\over{\Theta
(U(Q_0)-U(Q_{\infty})-e)s_{\infty}}}
\end{eqnarray*}
Since $e^{2u_i}=|h_{n,0}^{(1)}|$, the proposition follows. $\Box$

\section{Relation to immersed surface}\label{se:immersion}\setcounter{equation}{0}\setcounter{theorem}{0} \setcounter{definition}{0}
\setcounter{example}{0}
Let us now discuss the analytical description of a 2-dimensional
surface immersed in the ${\frak su}(N)$ algebra, associated with
the generalized Weierstrass formula for the immersion given below
by the expressions (\ref{eq:immerse}) and (\ref{eq:immersion}).

Suppose now one of the times, say $t_m$ in the hierarchy is
trivial, that is,
\begin{eqnarray*}
\p_{t_m}A_{\lambda}=\p_{t_m}B_{\lambda}=0
\end{eqnarray*}
then for the corresponding matrix $L_m$, we have the following
compatibility conditions
\begin{eqnarray}\label{eq:lax}
\p L_m&=&[L_m,A_{\lambda}]\nonumber\\
\overline{\p}L_m&=&[L_m,B_{\lambda}]
\end{eqnarray}
By applying the involution $\rho$ to equations (\ref{eq:lax}), we
get
\begin{eqnarray}\label{eq:star}
\p \left[\rho(L_m)\right]&=&[\rho(L_m),A_{\lambda}]\nonumber\\
\overline{\p} \left[\rho(L_m)\right]&=&[\rho(L_m),B_{\lambda}]
\end{eqnarray}
Therefore the matrix
\begin{eqnarray*}
X=i(L_m+\rho(L_m))\in{\frak su}(N+1)
\end{eqnarray*}
also satisfies the above equalities (\ref{eq:star}).

The spectral curve $\Sigma$ is the zero set of the determinant
\begin{eqnarray*}
\det (y-(L_m(\lambda)+\rho(L_m(\lambda)))
\end{eqnarray*}
It is a Riemann surface defined by a polynomial in $\lambda$ and
$y$ of the following form
\begin{eqnarray*}
y^{N+1}=a_m\lambda^m+a_{n-1}\lambda^{n-1}+\overline{a}_m\lambda^{-m}+P(y,\lambda)
\end{eqnarray*}
where $P(y,\lambda)$ is a polynomial in $\lambda$ and $y$. We see
that this curve has the form required by (\ref{eq:RS}). Therefore
we can apply the results of the previous section and express the
solution of the Toda lattice in terms of theta function on the
spectral curve.

If we set $|\lambda |=1$, we have the following compatibility
conditions
\begin{eqnarray}\label{eq:conserv}
\p(L_m+L_m^{\dag})=[L_m+L_m^{\dag},A_{\lambda}],\quad
\pbar(L_m+L_m^{\dag})=[L_m+L_m^{\dag},B_{\lambda}]
\end{eqnarray}
We can define the $N\times N$ matrix
\begin{eqnarray}\label{eq:immerse}
K(\lambda)&=&[L_m+L_m^{\dag},B_{\lambda}]
\end{eqnarray}
and its hermitian conjugate
\begin{eqnarray*}
K^{\dag}(\lambda)&=&[L_m+L_m^{\dag},A_{\lambda}]
\end{eqnarray*}
From (\ref{eq:conserv}) we obtain the following
\begin{eqnarray}\label{eq:conlaw}
\p K-\pbar K^{\dag}=0
\end{eqnarray}
Hence $\p K\in i{\frak su}(N)$ is a hermitian matrix.

As a consequence of the conservation law (\ref{eq:conlaw}), we
will show in the next section that there exists real-valued
functions $X_i(\xi,\overline{\xi})$ which can be identified with
the Weierstrass representation of surfaces immersed in the
multi-dimensional space $\mathbb{R}^{N^2-1}$.

By using the Drinfeld-Sokolov iteration process \cite{DS} introduced in
section \ref{se:DS}, the matrix
\begin{eqnarray*}
X=i(L_m+L_m^{\dag})
\end{eqnarray*}
and hence the immersion, can be computed.
\subsection{The generalized Weierstrass formula for immersion}
In order to study immersions defined by means of theta function
solutions of the Toda Lattice and to derive from this model the
moving frames and the corresponding Gauss-Weingarten and the
Gauss-Codozzi-Ricci equations, it is convenient to exploit the
Euclidean structure of the ${\frak su}(N)$ Lie algebra leading to
an identification
\begin{eqnarray*}
\mathbb{R}^{N^2-1}\cong {\frak su}(N)
\end{eqnarray*}
Consequently, we can introduce on ${\frak su}(N)$ an inner product
$(,):{\frak su}(N)\times{\frak su}(N)\rightarrow\mathbb{R}$ of
vectors in terms of matrices
\begin{eqnarray*}
(X,Y)=-{1\over 2}\tr(XY),\quad X,Y\in{\frak su}(N)
\end{eqnarray*}
Let us assume that the matrix $K$ is constructed from a solution
of the equation of motion written in the form of the conservation law
(\ref{eq:conlaw}). This conservation law implies that the
matrix-valued 1-form
\begin{eqnarray}\label{eq:form}
{\rm d}X=i(K^{\dag}\d\xi+K\d\overline{\xi})
\end{eqnarray}
is closed ($\d(\d X)=0$) and it takes values in the Lie algebra $\su
(N)$ of the anti-hermitian matrices. The real and the imaginary
parts of $\d X$ are anti-symmetric and symmetric, respectively
\begin{eqnarray*}
\d X&=&\d X^1+i\d X^2 \\
(\d X^1)^T&=&-\d X^1,\quad (\d X^2)^T=\d X^2
\end{eqnarray*}
where the 1-forms $\d X^1$ and $\d X^2$ take values in ${\frak
sl}(N,\mathbb{R})$.

From the closeness of $\d X$ it follows that the integral
\begin{eqnarray}\label{eq:immersion}
X(\xi,\overline{\xi})=i\int_{\gamma}(K^{\dag}\d\xi+K\d\overline{\xi})
\end{eqnarray}
locally depends only on the end points of the curve $\gamma$ in $\mathbb{C}$.

The integral defines the mapping
\begin{eqnarray*}
X:\Omega\in\mathbb{C}\rightarrow X(\xi,\overline{\xi})\in\su
(N)\cong\mathbb{R}^{N^2-1}
\end{eqnarray*}
which is called the generalized Weierstrass formula for immersion
(GWFI). The complex tangent vectors of this immersion, by virtue
of (\ref{eq:form}), are
\begin{eqnarray*}
\p X=iK^{\dag}, \quad \overline{\p}X=iK
\end{eqnarray*}
Hence a surface ${\cal F}$ associated with the Toda model by means
of the immersion (\ref{eq:immerse}) satisfies the following
relations
\begin{eqnarray*}
\p X&=&i[(L_m^{\dag}+L_m),B_{\lambda}] \\
\pbar X&=&i[(L_m^{\dag}+L_m),A_{\lambda}]
\end{eqnarray*}
The components of the metric on $\F$ induced by the Euclidean
structure in $\frak{su}(N)$ are given by
\begin{eqnarray*}
g_{11}&=&(\p X,\p X)={1\over
2}\tr(\p(L_m+L_m^{\dag})\p(L_m+L_m^{\dag}))\\
g_{12}&=&g_{21}=(\p X,\pbar X)=-{1\over
2}\tr(\p(L_m+L_m^{\dag})\pbar(L_m+L_m^{\dag}))\\
g_{22}&=&(\pbar X,\pbar X)={1\over
2}\tr(\pbar(L_m+L_m^{\dag})\pbar(L_m+L_m^{\dag}))
\end{eqnarray*}
The first fundamental form of the surface $\F$ takes the form
\begin{eqnarray}\label{eq:first}
I=g_{11}(\d\xi)^2+2g_{21}\d\xi\d\overline{\xi}+g_{22}(\d\overline{\xi})^2
\end{eqnarray}
and the second fundamental form is
\begin{eqnarray}\label{eq:second}
II=(\p^2X)^{\perp}\d\xi\d\xi+2(\p\pbar
X)^{\perp}\d\xi\d\overline{\xi}+(\pbar^2
X)^{\perp}\d\overline{\xi}\d\overline{\xi}
\end{eqnarray}
where $A^{\perp}$ is the normal component of $A$ to the surface
$\F$. The mean curvature vector is given by
\begin{eqnarray*}
H=\det G^{-1}(g_{22}(\p^2X)^{\perp}-2g_{12}(\p\pbar
X)^{\perp}+g_{11}(\pbar^2 X)^{\perp})
\end{eqnarray*}
where $G$ is the matrix formed by the metric $G_{ij}=g_{ij}$.

The Gaussian curvature is not neccessarily constant
\begin{eqnarray*}
K&=&(2\det G)^{-1}\Bigg(\p\left({1\over{\det G}}({g_{12}}\pbar\ln
g_{11}-\p g_{22})\right)\\
&+&\pbar \left({1\over {\det G}}(2\p g_{12}-\pbar
g_{11}-g_{12}\p\ln g_{11})\right)\Bigg)
\end{eqnarray*}
It is easy to check that
\begin{eqnarray*}
\tr(\p\pbar L\p L)&=&\tr([ L,\pbar A][L,A])\\
\tr(\p\pbar L\pbar L)&=&\tr([L,\p B][L,B])=\tr([ L,\pbar A][L,B])\\
\end{eqnarray*}
By equation (\ref{eq:lax}), we see that
\begin{eqnarray*}
\p X=[X,A_{\lambda}],\quad \overline{\p}X=[X,B_{\lambda}]
\end{eqnarray*}
Hence we have
\begin{eqnarray}\label{eq:conq}
k\p \tr(X^k)=\tr(X^{k-1}[X,A_{\lambda}]), \quad
=\tr([X^k,A_{\lambda}])=0
\end{eqnarray}
therefore the invariants of the surface is generated by
$\tr(X^i)$, $i=2,\ldots, m$. In particular, these surfaces are all
mapped into hyperspheres immersed in $\mathbb{R}^{N^2-1}$.

\section{The immersion in the ${\frak su}(N)$ case}\label{se:suN}\setcounter{equation}{0}\setcounter{theorem}{0} \setcounter{definition}{0}
\setcounter{example}{0}

In the ${\frak su}(N)$ case, the matrix $L_n(\lambda)$ can be
computed as follows. First we make the gauge transformation
\begin{eqnarray*}
\hat{A}_{\lambda}&=&
\mathbb{G}^{-1}A_{\lambda}\mathbb{G}+\mathbb{G}^{-1}\p
\mathbb{G}=\pmatrix{2\p u_0&0 &0&\ldots &0&\lambda \cr
             1&2\p u_1&0&\ldots &0&0\cr
             0&1&2\p u_2&\ldots &0&0&\cr
             \vdots &\vdots &\vdots &\ddots &\vdots &\vdots\cr
             0&0&0&\ldots &2\p u_{N-1}&0\cr
             0&0&0&\ldots &1&2\p u_N\cr}\\
\mathbb{G}&=&\diag\left(e^{u_0},e^{u_1},\ldots,e^{u_N}\right)
\end{eqnarray*}
We can then expand $\hat{A}_{\lambda}$ and $L_n(\lambda)$ in terms
of $\Lambda$
\begin{eqnarray*}
\hat{A}_{\lambda}&=&\pmatrix{2\p u_0&0 &0&\ldots &0&0 \cr
             0&2\p u_1&0&\ldots &0&0\cr
             0&0&2\p u_2&\ldots &0&0&\cr
             \vdots &\vdots &\vdots &\ddots &\vdots &\vdots\cr
             0&0&0&\ldots &2\p u_{N-1}&0\cr
             0&0&0&\ldots &0&2\p u_N\cr}+\Lambda\\
L_n(\lambda)&=&\sum_{i=0}^{3n+1}l_{3n+1-i}\Lambda^i
\end{eqnarray*}
where $l_i$ are diagonal matrices.

The compatibility condition (\ref{eq:lax}) becomes
\begin{eqnarray}\label{eq:temp}
\p L_n(\lambda)&=&[L_n(\lambda),\Lambda+\hat{A}_0]\nonumber\\
\hat{A}_0&=&2\diag\left(\p u_0,\ldots,\p u_N\right)
\end{eqnarray}
The expansion of (\ref{eq:temp}) in terms of $\Lambda$ then takes
the following form
\begin{eqnarray*}
\sum_{i=1}^{3n+1}\p l_{3n+1-i}\Lambda^i=\sum_{i=1}^{3n+1}
(l_{3n+1-i}\Lambda^{i+1}-\Lambda
l_{3n+1-i}\Lambda^i-\hat{A}_0l_{3n+1-i}\Lambda^i+l_{3n+1-i}\Lambda^i\hat{A}_0)
\end{eqnarray*}
Let $\sigma$ be the permutation such that
\begin{eqnarray*}
(\sigma(0),\sigma(1),\ldots,\sigma(N))=(N,0,1,\ldots,N-1)
\end{eqnarray*}
then we have
\begin{eqnarray*}
\Lambda\pmatrix{a_0&\cdots&\cdots&0\cr
0&a_1\cdots&0\cr\vdots&\vdots&\ddots&\vdots\cr0&\cdots&\cdots&a_N\cr}=
\pmatrix{a_{\sigma(0)}&\cdots&\cdots&0\cr
0&a_{\sigma(1)}\cdots&0\cr\vdots&\vdots&\ddots&\vdots\cr0&\cdots&\cdots&a_{\sigma(N)}\cr}\Lambda
\end{eqnarray*}
By using this permutation, we can express (\ref{eq:temp}) in terms
of the coefficients $l_i$.
\begin{eqnarray}\label{eq:comcoef}
l_{i+1}-\sigma(l_{i+1})=\p l_{i}-l_{i}(\sigma^{3n+1-i}(A_0)-A_0)
\end{eqnarray}
where the action of $\sigma$ on a diagonal matrix is defined as
follows
\begin{eqnarray*}
\sigma\diag\left(m_0,m_2,\ldots,m_N\right)=\diag(m_{\sigma(0)},m_{\sigma(1)},\ldots,m_{\sigma(N)})
\end{eqnarray*}
The equation (\ref{eq:comcoef}) determines the coefficients $l_i$
uniquely up to an addition of a scalar provided the right hand
side is traceless. The condition that the right hand side of
(\ref{eq:comcoef}) is traceless then fixes the scalar freedom in
$l_i$.

For example, in the case of ${\frak su}(3)$ and $n=1$, we have the
explicit form of the coefficients $l_i$
\begin{eqnarray}\label{eq:coefsu3}
l_0&=&I_d\nonumber\\
l_1&=&\pmatrix{2\p u_0&0&0\cr 0&2\p u_1&0\cr 0&0&2\p u_2\cr} \nonumber\\
l_2&=&\pmatrix{{2\over 3}(\p^2 u_0-\p^2u_1-c_2)&0&0\cr 0&{2\over
3}(\p^2 u_1-2\p^2u_2-c_2)&0\cr
0&0&{2\over 3}(2\p^2 u_2-\p^2u_0-c_2)\cr} \nonumber\\
c_i&=&(\p u_0)^i+(\p u_1)^i+(\p u_2)^i \nonumber\\
l_3&=&\diag(l_3^1,l_3^2,l_3^3)\\
l_3^1&=&{1\over 3}\left(-2\p^3u_1-\p c_2-4\p^2u_0(\p u_1-\p u_0)+4\p u_1c_2-{4\over 3}c_3\right) \nonumber\\
l_3^2&=&{1\over 3}\left(-2\p^3u_2-3\p c_2-4\p^2u_0(\p u_1-\p
u_0)-{4\over 3}c_3
+2\p(\p u_2-\p u_1)^2+4\p u_2c_2\right) \nonumber\\
l_3^3&=&{1\over 3}\left(-2\p^3u_0+\p c_2-4\p^2u_0(\p u_1-\p
u_0)-{4\over 3}c_3-2\p(\p u_0-\p u_1)^2-4\p u_0c_2\right) \nonumber\\
l_4&=&\diag(1_4^1,l_4^2,l_4^3) \nonumber\\
l_4^1&=&-{1\over 3}\left((\p l_3^1+2(\p u_0-\p u_2)l_3^1)-(\p
l_3^2+2(\p u_1-\p u_0)l_3^2)\right) \nonumber\\
l_4^2&=&{1\over 3}\left(-(\p l_3^1+2(\p u_0-\p u_2)l_3^1)+2(\p
l_3^2+2(\p u_1-\p u_0)l_3^2)\right) \nonumber\\
l_4^3&=&{1\over 3}\left(-2(\p l_3^1+2(\p u_0-\p u_2)l_3^1)+(\p
l_3^2+2(\p u_1-\p u_0)l_3^2)\right) \nonumber
\end{eqnarray}

\subsection{Moving frame in the ${\frak su}(N)$ case}

In this section we will construct a moving frame for the immersed
surface ${\cal F}$ in the ${\frak su}(N)$ case, which will then be
used to compute the Gauss-Weingarten equations.

To begin with, let the action of the involution $\rho$
(\ref{eq:invol}) on a matrix-valued function $\Phi(P)$ on the
spectral curve be defined as follows
\begin{eqnarray}\label{eq:phoaction}
\rho(\Phi(P))=\Phi^{\dag}(\rho(P))
\end{eqnarray}
We then have the following
\begin{lemma}\label{le:inverse}
Let $F_1$ be the matrix in (\ref{eq:F1}), and $\Phi=\mathbb{G}
F_1$, where
\begin{eqnarray*}
\mathbb{G}=\diag\left(e^{u_0},e^{u_1},\ldots,e^{u_N}\right)
\end{eqnarray*}
then
\begin{eqnarray}\label{eq:inverse}
\rho(\Phi)\Phi&=&\diag\left(d_1
,\ldots,d_{N+1}\right)\nonumber\\
&=&\mathbb{D}\nonumber\\
d_j&=&\sum_{i=0}^{N}e^{u_i+\overline{u_i}}f_i(P_j)\overline{f}_i(\rho(P_j))
\end{eqnarray}
\end{lemma}

Proof. Let $\Phi=\mathbb{G}F_1$, $H_1=\rho(\Phi)$ and
$H_2=(\Phi)^{-1}$. First note that both $H_1$ and $H_2$ are
matrices that has the following asymptotic behavior
\begin{eqnarray}\label{eq:asym}
H_j&=&\exp D_0T_0^j\lambda^S,\quad \lambda\rightarrow 0\nonumber\\
H_j&=&\exp D_{\infty}T_{\infty}^j\lambda^S,\quad\lambda\rightarrow\infty\nonumber\\
D_0&=&-\xi\lambda^{1\over
N+1}\diag\left(1,\omega,\omega^2,\ldots,\omega^N\right)\\
D_{\infty}&=&\overline{\xi}\lambda^{-{1\over
N+1}}\diag\left(1,\omega,\omega^2,\ldots,\omega^N\right)\nonumber\\
S&=&\diag\left(0,{1\over {N+1}},\ldots,{N\over {N+1}}\right)\nonumber\\
\omega&=&e^{{2\pi i}\over{N+1}}\nonumber
\end{eqnarray}
where $T_{0}^i$ and $T_{\infty}^i$ are power series in
$\lambda^{1\over {N+1}}$ and $\lambda^{-{1\over {N+1}}}$ that are
invertible at $0$ and $\infty$ respectively.

By applying the involution $\rho$ to the equations
\begin{eqnarray*}
\p\Phi+A_{\lambda}\Phi=0,\quad \overline{\p}\Phi+B_{\lambda}\Phi=0
\end{eqnarray*}
we obtain
\begin{eqnarray}\label{eq:dual}
\overline{\p}H_1-H_1B_{\lambda}=0,\quad \p H_1-H_1A_{\lambda}=0
\end{eqnarray}
Similarly, by differentiating $\Phi^{-1}$, we see that $H_2$
satisfies the same equations
\begin{eqnarray*}
\overline{\p}H_2-H_2B_{\lambda}=0,\quad \p H_2-H_2A_{\lambda}=0
\end{eqnarray*}
Since any solution to the equations (\ref{eq:dual}) with the
asymptotic behavior (\ref{eq:asym}) are determined uniquely up to
the multiplication of a diagonal matrix $\mathbb{D}$ on the left, we
have
\begin{eqnarray*}
H_1=\mathbb{D} H_2
\end{eqnarray*}
for some diagonal matrix $\mathbb{D}$ constant in $\xi$ and
$\overline{\xi}$. Therefore we have
\begin{eqnarray*}
H_1\Phi =\rho(F_1)\overline{\mathbb{G}}\mathbb{G}F_1=\mathbb{D}
\end{eqnarray*}
By computing the diagonal entries in
$\rho(F_1)\overline{\mathbb{G}}\mathbb{G}F_1$, the lemma is
proven. $\Box$

We can now construct the moving frame for the surface ${\cal F}$
immersed in ${\frak su}(N+1)$. First note that the function
$T(P)=\Phi(P)\exp(-D_0-D_{\infty})$, where $\Phi$, $D_0$ and
$D_{\infty}$ are defined in (\ref{eq:asym}), satisfies the
following
\begin{eqnarray*}
\p TT^{-1}+A_{\lambda}=TD_AT^{-1},\quad \overline{\p}
TT^{-1}+B_{\lambda}=TD_BT^{-1}
\end{eqnarray*}
where $D_A$ and $D_B$ are diagonal matrices. We see that $T$ and
the dressing matrix $T_{\infty}$ in (\ref{eq:dressing}) are
related by
\begin{eqnarray*}
T=T_{\infty}\Psi
\end{eqnarray*}
where $\Psi$ is the matrix that diagonalizes $\Lambda$.

Since $L_m$ in (\ref{eq:lax}) are defined by (\ref{eq:Li})
\begin{eqnarray*}
L_i=(T_{\infty}\Lambda^iT_{\infty}^{-1})_+=(TDT^{-1})+\p_m
TT^{-1}=(TDT^{-1})
\end{eqnarray*}
for some diagonal matrix $D$ that is constant in $\xi$ and
$\overline{\xi}$. The second equality follows as $\p_mT=0$.
Therefore $T$, and hence $\Phi$, diagonalizes $L_m$, and hence
$X$.

We can now take the conjugation of $A_{\lambda}$, $B_{\lambda}$
and $X$ by the matrix $\Phi$ to obtain the following
\begin{eqnarray}\label{eq:adjoint}
X&=&\Phi Y\Phi^{-1}\nonumber\\
Y&=&\diag\left(y_1,\ldots,y_{N+1}\right) \nonumber\\
A_{\lambda}&=&\Phi\Xi_1\Phi^{-1},\quad
B_{\lambda}=\Phi\Xi_2\Phi^{-1}\nonumber\\
\left(\Xi_1\right)_{kl}&=&-\left(\sum_{i=0}^{N}\phi_i(P_k)\overline{\phi}_i(\rho(P_k))\right)^{-1}
\sum_{i=0}^{N}\overline{\phi}_i(\rho(P_k))\p\phi_i(P_l)\\
\left(\Xi_2\right)_{kl}&=&-\left(\sum_{i=0}^{N}\phi_i(P_k)\overline{\phi}_i(\rho(P_k))\right)^{-1}
\sum_{i=0}^{N}\overline{\phi}_i(\rho(P_k))\overline{\p}\phi_i(P_l)\nonumber\\
\phi_i(P)&=&e^u_if_i(P)\nonumber
\end{eqnarray}
the entries $y_i$ of $Y$ are the different branches of the spectrum $y$
in the spectral curve (\ref{eq:RS}). Note that they do not depend on $\xi$
and $\overline{\xi}$.

The tangent vectors $\p X$ and $\overline{\p}X$ can be represented
as
\begin{eqnarray}\label{eq:tangent}
\p X&=&[X,A_{\lambda}]=\Phi \p X^0\Phi^{-1},\quad \overline{\p}
X=[X,B_{\lambda}]=\Phi \overline{\p} X^0\Phi^{-1}\nonumber\\
\left(\p X^0\right)_{kl}&=&(y_k-y_l)\left(\Xi_1\right)_{kl}\\
\left(\overline{\p}
X^0\right)_{kl}&=&(y_k-y_l)\left(\Xi_2\right)_{kl}\nonumber
\end{eqnarray}
From this the metric and the first fundamental form are
\begin{eqnarray}\label{eq:1stform}
g_{11}&=&-{1\over 2}\tr(\p X\p X)=
-{1\over 2}\sum_{k,l=1}^{N+1}(y_k-y_l)^2(\Xi_1)_{kl}(\Xi_1)_{lk}\nonumber\\
g_{12}&=&g_{21}=-{1\over 2}\tr(\p X\pbar
X)=-{1\over 2}\sum_{k,l=1}^{N+1}(y_k-y_l)^2(\Xi_1)_{kl}(\Xi_2)_{lk}\\
g_{22}&=&-{1\over 2}\tr(\pbar X\pbar X)=
-{1\over 2}\sum_{k,l=1}^{N+1}(y_k-y_l)^2(\Xi_2)_{kl}(\Xi_2)_{lk}\nonumber\\
I&=&g_{11}\d\xi\d\xi+2g_{12}\d\xi\d\overline{\xi}+g_{22}\d\overline{\xi}\d\overline{\xi}\nonumber
\end{eqnarray}
Let $E_{ij}$ be the $N+1\times N+1$ matrix that is 1 in its
$ij^{th}$ entry and zero elsewhere. The following forms a basis of
the ${\frak su}(N+1)$ algebra
\begin{eqnarray*}\label{eq:basis}
s(k,l)&=&{{(k-2)(k-1)}\over 2}+l\nonumber\\
v_{s(k,l)}^0&=&E_{kl}-E_{lk}, \quad l<k,\quad k=2,\ldots, N+1 \nonumber\\
u_{s(k,l)}^0&=&i(E_{kl}+E_{lk}), \quad l<k,\quad k=2,\ldots, N+1\nonumber \\
d_{k}^0&=&i(E_{11}-E_{kk}),\quad k=2,\ldots, N+1\nonumber\\
v_{s(k,l)}&=&\Phi v_{s(k,l)}^0\Phi^{-1}\\
u_{s(k,l)}&=&\Phi u_{s(k,l)}^0\Phi^{-1}\nonumber\\
d_{k}&=&\Phi d_{k}^0\Phi^{-1}\nonumber
\end{eqnarray*}
We could now construct an orthonormal basis to the surface in
terms of these basis vectors.

\begin{theorem}
The following forms an orthogonal basis of the normal vectors to
the surface ${\cal F}$.
\begin{eqnarray}\label{eq:frame}
U_{l}&=&\sum_{i=1}^{l}C_i^lu_i,\quad l=3,\ldots, {{N(N+1)}\over 2}\nonumber\\
V_{l}&=&\sum_{i=1}^{l}D_i^lv_i,\quad l=3,\ldots, {{N(N+1)}\over 2}\nonumber\\
W_0&=&X\\
W_l&=&\sum_{i=2}^{l}K_i^ld_i\nonumber
\end{eqnarray}
where the constants $C_i^l$, $D_i^l$ and $K_i^l$ are given by the
following
\begin{eqnarray*}
C_i^l&=&0, i>l\nonumber\\
C_i^l&=&\det\pmatrix{J_1&J_2&\ldots
&J_{i-1}&-J_{l}&J_{i+1}&\ldots&J_{l-1}\cr
                 -\overline{J}_1&-\overline{J}_2&\ldots
                 &-\overline{J}_{i-1}&\overline{J}_{l}&-\overline{J}_{i+1}&\ldots&-\overline{J}_{l-1}\cr
                 C_1^{l-1}&C_2^{l-1}&\ldots &C_{i-1}^{l-1}&0&C_{i+1}^{l-1}&\ldots&C_{l-1}^{l-1}\cr
                 \vdots&\vdots&\vdots&\ddots&\vdots&\vdots&\vdots&\vdots\cr
                 C_1^{3}&C_2^3&\ldots
                 &C_{i-1}^{3}&0&C_{i+1}^{3}&\ldots&0\cr}\nonumber\\
C_l^l&=&\det\pmatrix{J_1&J_2&\ldots&J_{l-1}\cr
                 -\overline{J}_1&-\overline{J}_2&\ldots&-\overline{J}_{l-1}\cr
                 C_1^{l-1}&C_2^{l-1}&\ldots&C_{l-1}^{l-1}\cr
                 \vdots&\vdots&\ddots&\vdots\cr
                 C_1^{3}&C_2^3&\ldots&0\cr}\\
D_i^l&=&\det\pmatrix{Q_1&Q_2&\ldots
&Q_{i-1}&-Q_{l}&Q_{i+1}&\ldots&Q_{l-1}\cr
                 -\overline{Q}_1&-\overline{Q}_2&\ldots
                 &-\overline{Q}_{i-1}&\overline{Q}_{l}&-\overline{Q}_{i+1}&\ldots&-\overline{Q}_{l-1}\cr
                 D_1^{l-1}&D_2^{l-1}&\ldots &D_{i-1}^{l-1}&0&D_{i+1}^{l-1}&\ldots&D_{l-1}^{l-1}\cr
                 \vdots&\vdots&\vdots&\ddots&\vdots&\vdots&\vdots&\vdots\cr
                 D_1^{3}&D_2^3&\ldots&D_{i-1}^{3}&0&D_{i+1}^{3}&\ldots&0\cr}\nonumber\\
D_l^l&=&\det\pmatrix{Q_1&Q_2&\ldots &Q_{l-1}\cr
                 -\overline{Q}_1&-\overline{Q}_2&\ldots&-\overline{Q}_{l-1}\cr
                 D_1^{l-1}&D_2^{l-1}&\ldots&D_{l-1}^{l-1}\cr
                 \vdots&\vdots&\ddots&\vdots\cr
                 D_1^{3}&D_2^3&\ldots&0\cr}\nonumber\\
J_{s(k,l)}&=&-{1\over 2}\tr(u_{s(k,l)}\p X)=-{i\over 2}((\p X^0)_{lk}+(\p X^0)_{kl})\nonumber\\
Q_{s(k,l)}&=&-{1\over 2}\tr(v_{s(k,l)}\p X)=-{1\over 2}((\p X^0)_{lk}-(\p X^0)_{kl})\nonumber\\
K_k^l&=&\det\pmatrix{R_1&R_2&\ldots
&R_{k-1}&-R_{l}&R_{k+1}&\ldots&R_{l-1}\cr
Z_1&Z_2&\ldots&Z_{k-1}&-1&Z_{k+1}&\ldots&Z_{l-1}\cr
\vdots&\vdots&\ddots&\vdots&\vdots&\vdots&\vdots&\vdots\cr
-2K_2^{2}&0&\ldots&0&-1& 0&\ldots&0\cr} \nonumber\\
Z_1&=&-2K_2^{l-1}+\sum_{j=3}^{l-1}K_j^{l-1}\nonumber\\
Z_k&=&-2K_{k-1}^{l-1}+K_2^{l-1}\nonumber\\
 R_k&=&\tr(d_kX)=i(y_1-y_k)
\end{eqnarray*}

\end{theorem}
The coefficients in (\ref{eq:frame}) can be computed by solving
systems of linear equations. To compute the first coefficient
$C_1^{l}$, we need to solve
\begin{eqnarray*}
a_1\tr(\p Xv_1)+a_2\tr(\p X v_2)+\tr(\p X v_3)&=&0\\
a_1\tr(\pbar Xv_1)+a_2\tr(\pbar X v_2)+\tr(\pbar X v_3)&=&0
\end{eqnarray*}
for $a_1$ and $a_2$. The solution is given by
\begin{eqnarray*}
a_1&=&-{\det\pmatrix{J_3&J_2\cr
\overline{J}_3&\overline{J}_2\cr}\over\det\pmatrix{J_1&J_2\cr
\overline{J}_1&\overline{J}_2\cr}} \\
a_2&=&-{\det\pmatrix{J_1&J_3\cr
\overline{J}_1&\overline{J}_3\cr}\over\det\pmatrix{J_1&J_2\cr
\overline{J}_1&\overline{J}_2\cr} }
\end{eqnarray*}
Then, by multiplying the vector
\begin{eqnarray*}
a_1v_1+a_2v_2+v_3
\end{eqnarray*}
by the common dominator of $a_1$ and $a_2$, one sees that
\begin{eqnarray*}
V_1=-\det\pmatrix{J_3&J_2\cr
\overline{J}_3&\overline{J}_2\cr}v_1-\det\pmatrix{J_1&J_3\cr
\overline{J}_1&\overline{J}_3\cr}v_2+\det\pmatrix{J_1&J_2\cr
\overline{J}_1&\overline{J}_2\cr}v_3
\end{eqnarray*}
is a normal vector to the surface ${\cal F}$. To compute the other
coefficients, one needs to solve linear systems of equations
\begin{eqnarray*}
a_1\tr(\p Xv_1)+\ldots+a_{l+1}\tr(\p Xv_{l-1})+\tr(\p X v_l)&=&0\\
a_1\tr(\pbar Xv_1)+\ldots+a_{l+1}\tr(\pbar X v_{l-1})+\tr(\pbar X
v_l)&=&0\\
a_1\tr(V_1v_1)+\ldots+a_{l+1}\tr(V_1v_{l-1})+\tr(V_1v_l)&=&0\\
\ldots\\
a_1\tr(V_{l-1}v_1)+\ldots+a_{l+1}\tr(V_{l-1}v_{l-1})+\tr(V_{l-1}v_l)&=&0
\end{eqnarray*}
and similar equations for the vectors $U_i$ and $W_i$.

\subsection{The Weingarten equation in a non-orthonormal basis}

We can construct another basis of the normal vectors as follows.
Let $a_{i1}$ and $a_{i2}$ be the solution of the following linear
equations
\begin{eqnarray*}
a_{i1}\tr(u_1\p X)+a_{i2}\tr(v_1\p X)+\tr(\p X v_i)&=&0\\
a_{i1}\tr(u_1\pbar X)+a_{i2}\tr(v_1\pbar X)+\tr(\pbar X v_i)&=&0
\end{eqnarray*}
Then $a_{i1}$ and $a_{i2}$ are given by the following
\begin{eqnarray*}
a_{i1}&=&-{\det\pmatrix{Q_i&Q_1\cr
\overline{Q}_i&\overline{Q}_1\cr}\over\det\pmatrix{J_1&Q_1\cr
\overline{J}_1&\overline{Q}_1\cr}} \\
a_{i2}&=&-{\det\pmatrix{J_1&Q_i\cr
\overline{J}_1&\overline{Q}_i\cr}\over\det\pmatrix{J_1&Q_1\cr
\overline{J}_1&\overline{Q}_1\cr} }
\end{eqnarray*}
By similar argument as before, we see that the vectors
\begin{eqnarray}\label{eq:normalvect}
n_{j}^v&=&\left(\gamma_{1,1}v_{j+1}-\gamma_{1,j+1}v_1-\kappa_{j+1,1}u_1\right),\quad j=1,\ldots {{N(N-1)}\over 2}\nonumber\\
n_{j}^u&=&\left(\gamma_{1,1}u_{j+1}-\beta_{1,j+1}v_1-\gamma_{j+1,1}u_1\right),\quad j=1,\ldots {{N(N-1)}\over 2}\nonumber\\
\beta_{k,l}&=&i\det\pmatrix{J_k&J_l\cr
\overline{J}_k&\overline{J}_l\cr}\\
\gamma_{k,l}&=&i\det\pmatrix{J_k&Q_l\cr
\overline{J}_k&\overline{Q}_l\cr}\nonumber\\
\kappa_{k,l}&=&i\det\pmatrix{Q_k&Q_l\cr
\overline{Q}_k&\overline{Q}_l\cr}\nonumber\\
 n_{j}^d&=&i\Phi(E_{11}-E_{j+1,j+1})\Phi^{-1},\quad
j=1,\ldots,
N-1\nonumber\\
n_{N}^d&=&iy_{N+1}^{-1}X\nonumber
\end{eqnarray}
form a basis of normal vectors to the surface ${\cal F}$. Although
this basis is not orthogonal, they are a lot simpler than the
orthonormal basis (\ref{eq:frame}) because the coefficients only
involves determinants of $2\times 2$ matrices. Note that $n_j^v$,
$n_j^u$ and $n_j^d$ all belong to ${\frak su}(N+1)$.

We will now calculate the Gauss-Weingarten equation with this
basis of normal vectors. By differentiating the vector $n_{j}^v$
with respect to $\xi$ we see that
\begin{eqnarray*}
\p
n_{j}^v&=&[\p\Phi\Phi^{-1},n_{j}^v]+\left(\p\gamma_{1,1}v_{j+1}-\p\gamma_{1,j+1}v_1-\p\kappa_{j+1,1}u_1\right)\\
&=&[n_{j}^v,A_{\lambda}]+\left(\p\gamma_{1,1}v_{j+1}-\p\gamma_{1,j+1}v_1-\p\kappa_{j+1,1}u_1\right)
\end{eqnarray*}
since $\p\Phi\Phi^{-1}=-A_{\lambda}$. By solving similar linear
equations, we see that the tangent component of the $\p n_{j}^v$
is given by
\begin{eqnarray}\label{eq:tang}
(\p n_{j}^v)_T&=&\nu_{j,1}\p X+\nu_{j,2}\pbar X\nonumber\\
\nu_{j1}&=&\det G^{-1}\det\pmatrix{b_{j1}&g_{11}\cr
b_{j2}&g_{12}\cr}\nonumber\\
\nu_{j2}&=&\det G^{-1}\det\pmatrix{g_{12}&b_{j1}\cr
\overline{g}_{11}&b_{j2}\cr}\\
b_{j1}&=&\p\gamma_{1,1}Q_{j+1}-\p\gamma_{1,j+1}Q_1-\p\kappa_{j+1,1}J_1-{1\over
2}\tr([n_{j}^v,A_{\lambda}]\p
X)\nonumber\\
b_{j2}&=&\p\gamma_{1,1}\overline{Q}_{j+1}-\p\gamma_{1,j+1}\overline{Q}_1-\p\kappa_{j+1,1}\overline{J}_1-{1\over
2}\tr([n_{j}^v,A_{\lambda}]\pbar X)\nonumber
\end{eqnarray}
Let $\eta$ be a normal vector to the surface ${\cal F}$, then
$\eta$ can be written in terms of the basis (\ref{eq:normalvect})
as follows
\begin{eqnarray*}
\eta&=&\sum_{j=1}^{{N(N-1)}\over
2}\eta_j^vn_j^v+\sum_{j=1}^{{N(N-1)}\over
2}\eta_j^un_j^u+\sum_{j=1}^{N}\eta_j^dn_j^d\\
&=&\sum_{j=1}^{{N(N-1)}\over
2}\gamma_{1,1}\eta_j^vv_{j+1}+\sum_{j=1}^{{N(N-1)}\over
2}\gamma_{1,1}\eta_j^uu_{j+1}+c_0^uu_1+c_0^vv_1\\
&+&\Phi\left(\sum_{j=1}^{N}i\eta_j^dE_{j+1,j+1}+c_1^dE_{11}\right)\Phi^{-1}
\end{eqnarray*}
Therefore the coefficients in this decomposition are given by
\begin{eqnarray}\label{eq:normalcoef}
\eta_i^v&=&-{1\over {2\gamma_{1,1}}}\tr(\eta v_{i+1}),\quad
\eta_i^u=-{1\over {2\gamma_{1,1}}}\tr(\eta u_{i+1}),\nonumber\\
\eta_j^d&=&-i\tr(\eta \Phi E_{j+1,j+1}\Phi^{-1})
\end{eqnarray}
By subtracting the tangent component (\ref{eq:tang}) from the
vector $\p n_j^v$, we obtain the normal component of $\p n_j^v$
\begin{eqnarray*}
(\p n_j^v)_{\perp}=
[n_{j}^v,A_{\lambda}]+\p\gamma_{1,1}v_{j+1}-\p\gamma_{1,j+1}v_1-\p\kappa_{j+1,1}u_1
-\nu_{j1}\p X-\nu_{j2}\pbar X
\end{eqnarray*}
By using (\ref{eq:normalcoef}), the coefficients of this vector
are given by the following
\begin{eqnarray}\label{eq:njvcoef}
(\p n_j^v)_k^v&=&\gamma_{1,1}^{-1}\left(-{1\over
2}\tr([n_{j}^v,A_{\lambda}]v_{k+1})+\p\gamma_{1,1}\delta_{kj}-\nu_{j1}Q_{k+1}-\nu_{j2}\overline{Q}_{k+1}\right)\nonumber\\
(\p n_j^v)_k^u&=&\gamma_{1,1}^{-1}\left(-{1\over
2}\tr([n_{j}^v,A_{\lambda}]u_{k+1})-\nu_{j1}J_{k+1}-\nu_{j2}\overline{J}_{k+1}\right)\\
(\p n_j^v)_k^d&=&-i\tr([n_{j}^v,A_{\lambda}]\Phi
E_{k+1,k+1}\Phi^{-1})\nonumber
\end{eqnarray}
We can perform similar computations with derivatives of other
normal vectors to obtain the Gauss-Weingarten equation.
\begin{theorem}\label{thm:GaussW}
The Gauss-Weingarten equation in terms of the basis
\begin{eqnarray*}
\eta=(\p X,\pbar X,n_j^v,n_j^u,n_j^d)^{T}
\end{eqnarray*} defined in (\ref{eq:normalvect}) is given by
\begin{eqnarray*}
\p^2X&=&\alpha_{1,1}\p X+\alpha_{1,2}\pbar
X+\sum_{k=1}^{{N(N-1)}\over
2}\left((\p^2X)_k^vn_k^v+(\p^2X)_k^un_k^u\right)\nonumber\\
&+&\sum_{k=1}^{N}(\p^2X)_k^dn_k^d\nonumber\\
\p\pbar X&=&\alpha_{2,1}\p X+\alpha_{2,2}\pbar
X+\sum_{k=1}^{{N(N-1)}\over 2}\left((\p\pbar
X)_k^vn_k^v\nonumber+(\p\pbar X)_k^un_k^u\right)\nonumber\\
&+&\sum_{k=1}^{N}(\p\pbar X)_k^dn_k^d\nonumber\\
\p n_j^v&=&\nu_{j,1}\p X+\nu_{j,2}\pbar
X+\sum_{k=1}^{{N(N-1)}\over 2}\left((\p
n_j^v)_k^vn_k^v+(\p n_j^v)_k^un_k^u\right)\nonumber\\
&+&\sum_{k=1}^{N}(\p n_j^v)_k^dn_k^d,\quad j=1,\ldots,{{N(N-1)}\over 2}\\
\p n_j^u&=&\mu_{j,1}\p X+\mu_{j,2}\pbar
X+\sum_{k=1}^{{N(N-1)}\over 2}\left((\p
n_j^u)_k^vn_k^v+(\p n_j^u)_k^un_k^u\right)\nonumber\\
&+&\sum_{k=1}^{N}(\p n_j^u)_k^dn_k^d,\quad j=1,\ldots,{{N(N-1)}\over 2}\nonumber\\
\p n_j^d&=&\chi_{j,1}\p X+\chi_{j,2}\pbar
X+\sum_{k=1}^{{N(N-1)}\over 2}\left((\p n_j^d)_k^vn_k^v+(\p
n_j^d)_k^un_k^u\right)\nonumber\\
&+&\sum_{k=1}^{N}(\p n_j^d)_k^dn_k^d,\quad j=1,\ldots,
N-1\nonumber\\
\p n_j^{N}&=&iy_{N+1}^{-1}\p X\nonumber
\end{eqnarray*}
where the coefficients are
\begin{eqnarray*}
\alpha_{1,1}&=&(\det G)^{-1}({1\over 2}\p g_{11}\overline{g}_{11}-\p g_{12}g_{12}+{1\over 2}\pbar g_{11}g_{12})\nonumber\\
\alpha_{1,2}&=&(\det G)^{-1}(g_{11}\p g_{12}-{1\over 2}g_{11}\pbar
g_{11}-{1\over 2}g_{12}\p g_{11})\nonumber \\
\alpha_{2,1}&=&(2\det G)^{-1}(\pbar g_{11}\overline{g}_{11}-\p \overline{g}_{11}g_{12})\nonumber\\
\alpha_{2,2}&=&(2\det G)^{-1}(g_{11}\p\overline{g}_{11}-\pbar
{g}_{11}g_{12})\nonumber\\
\nu_{j1}&=&\det G^{-1}\det\pmatrix{b_{j1}&g_{11}\cr
b_{j2}&g_{12}\cr}\nonumber\\
\nu_{j2}&=&\det G^{-1}\det\pmatrix{g_{12}&b_{j1}\cr
\overline{g}_{11}&b_{j2}\cr}\nonumber\\
b_{j1}&=&\p\gamma_{1,1}Q_{j+1}-\p\gamma_{1,j+1}Q_1-\p\kappa_{j+1,1}J_1-{1\over
2}\tr([n_{j}^v,A_{\lambda}]\p
X)\nonumber\\
b_{j2}&=&\p\gamma_{1,1}\overline{Q}_{j+1}-\p\gamma_{1,j+1}\overline{Q}_1-\p\kappa_{j+1,1}\overline{J}_1-{1\over
2}\tr([n_{j}^v,A_{\lambda}]\pbar X)\nonumber\\
\mu_{j1}&=&\det G^{-1}\det\pmatrix{a_{j1}&g_{11}\cr
a_{j2}&g_{12}\cr}\nonumber\\
\mu_{j2}&=&\det G^{-1}\det\pmatrix{g_{12}&a_{j1}\cr
\overline{g}_{11}&a_{j2}\cr}\nonumber\\
a_{j1}&=&\p\gamma_{1,1}J_{j+1}-\p\beta_{1,j+1}Q_1-\p\gamma_{j+1,1}J_1-{1\over
2}\tr([n_{j}^u,A_{\lambda}]\p
X)\nonumber\\
a_{j2}&=&\p\gamma_{1,1}\overline{J}_{j+1}-\p\beta_{1,j+1}\overline{Q}_1-\p\gamma_{j+1,1}\overline{J}_1-{1\over
2}\tr([n_{j}^v,A_{\lambda}]\pbar X)\nonumber\\
\chi_{j1}&=&\det G^{-1}\det\pmatrix{-{1\over
2}\tr([n_j^d,A_{\lambda}]\p X)&g_{11}\cr -{1\over
2}\tr([n_j^d,A_{\lambda}]\pbar X)&g_{12}\cr}\nonumber\\
\chi_{j2}&=&\det G^{-1}\det\pmatrix{g_{12}&-{1\over
2}\tr([n_j^d,A_{\lambda}]\p X)\cr \overline{g}_{11}&-{1\over
2}\tr([n_j^d,A_{\lambda}]\pbar X)\cr}\nonumber\\
(\p^2X)_k^v&=&\gamma_{1,1}^{-1}\left(\p Q_{j+1}+{1\over 2}\tr(\p
X[v_{j+1},A_{\lambda}])-\alpha_{1,1}Q_{j+1}-\alpha_{1,2}\overline{Q}_{j+1}\right)\nonumber\\
(\p^2X)_k^u&=&\gamma_{1,1}^{-1}\left(\p J_{j+1}+{1\over 2}\tr(\p
X[u_{j+1},A_{\lambda}])-\alpha_{1,1}J_{j+1}-\alpha_{1,2}\overline{J}_{j+1}\right)\nonumber\\
(\p^2X)_k^d&=&-i\tr(\p^2X\Phi
E_{k+1,k+1}\Phi^{-1})\nonumber\\
(\p\pbar X)_k^v&=&\gamma_{1,1}^{-1}\left(\p
\overline{Q}_{j+1}+{1\over 2}\tr(\pbar
X[v_{j+1},A_{\lambda}])-\alpha_{2,1}Q_{j+1}-\alpha_{2,2}\overline{Q}_{j+1}\right)\nonumber\\
(\p\pbar X)_k^u&=&\gamma_{1,1}^{-1}\left(\p
\overline{J}_{j+1}+{1\over 2}\tr(\pbar
X[u_{j+1},A_{\lambda}])-\alpha_{2,1}J_{j+1}-\alpha_{2,2}\overline{J}_{j+1}\right)\nonumber\\
(\p\pbar X)_k^d&=&-i\tr(\p\pbar X\Phi
E_{k+1,k+1}\Phi^{-1})\nonumber\\
(\p n_j^v)_k^v&=&\gamma_{1,1}^{-1}\left(-{1\over
2}\tr([n_{j}^v,A_{\lambda}]v_{k+1})+\p\gamma_{1,1}\delta_{kj}-\nu_{j1}Q_{k+1}-\nu_{j2}\overline{Q}_{k+1}\right)\nonumber\\
(\p n_j^v)_k^u&=&\gamma_{1,1}^{-1}\left(-{1\over
2}\tr([n_{j}^v,A_{\lambda}]u_{k+1})-\nu_{j1}J_{k+1}-\nu_{j2}\overline{J}_{k+1}\right)\\
(\p n_j^v)_k^d&=&-i\tr([n_{j}^v,A_{\lambda}]\Phi
E_{k+1,k+1}\Phi^{-1})\nonumber \\
(\p n_j^u)_k^v&=&\gamma_{1,1}^{-1}\left(-{1\over
2}\tr([n_{j}^u,A_{\lambda}]v_{k+1})-\mu_{j1}Q_{k+1}-\mu_{j2}\overline{Q}_{k+1}\right)\nonumber\\
(\p n_j^u)_k^u&=&\gamma_{1,1}^{-1}\left(-{1\over
2}\tr([n_{j}^u,A_{\lambda}]u_{k+1})+\p\gamma_{1,1}\delta_{jk}-\mu_{j1}J_{k+1}-\mu_{j2}\overline{J}_{k+1}\right)\nonumber\\
(\p n_j^u)_k^d&=&-i\tr([n_{j}^u,A_{\lambda}]\Phi
E_{k+1,k+1}\Phi^{-1})\nonumber\\
(\p n_j^d)_k^v&=&\gamma_{1,1}^{-1}\left(-{1\over
2}\tr([n_{j}^d,A_{\lambda}]v_{k+1})-\chi_{j1}Q_{k+1}-\chi_{j2}\overline{Q}_{k+1}\right)\nonumber\\
(\p n_j^d)_k^u&=&\gamma_{1,1}^{-1}\left(-{1\over
2}\tr([n_{j}^d,A_{\lambda}]u_{k+1})-\chi_{j1}J_{k+1}-\chi_{j2}\overline{J}_{k+1}\right)\nonumber\\
(\p n_j^d)_k^d&=&-i\tr([n_{j}^d,A_{\lambda}]\Phi
E_{k+1,k+1}\Phi^{-1})\nonumber
\end{eqnarray*}
The other equations can be obtained by taking hermitian
conjugation and by using the fact that $X$, $n_j^v$, $n_j^u$ and
$n_j^d$ all belong to ${\frak su}(N+1)$.
\end{theorem}
The second fundamental form (\ref{eq:second}) can now be read off from the
Gauss-Weingarten equation
\begin{eqnarray*}
(\p^2X)_{\perp}&=&\sum_{k=1}^{{N(N-1)}\over
2}\left((\p^2X)_k^vn_k^v+(\p^2X)_k^un_k^u\right)\nonumber\\
&+&\sum_{k=1}^{N}(\p^2X)_k^dn_k^d\nonumber\\
(\p\pbar X)_{\perp}&=&\sum_{k=1}^{{N(N-1)}\over 2}\left((\p\pbar
X)_k^vn_k^v\nonumber+(\p\pbar X)_k^un_k^u\right)\nonumber\\
&+&\sum_{k=1}^{N}(\p\pbar X)_k^dn_k^d\nonumber\\
(\p^2X)_k^v&=&\gamma_{1,1}^{-1}\left(\p Q_{j+1}+{1\over 2}\tr(\p
X[v_{j+1},A_{\lambda}])-\alpha_{1,1}Q_{j+1}-\alpha_{1,2}\overline{Q}_{j+1}\right)\nonumber\\
(\p^2X)_k^u&=&\gamma_{1,1}^{-1}\left(\p J_{j+1}+{1\over 2}\tr(\p
X[u_{j+1},A_{\lambda}])-\alpha_{1,1}J_{j+1}-\alpha_{1,2}\overline{J}_{j+1}\right)\nonumber\\
(\p^2X)_k^d&=&-i\tr(\p^2X\Phi
E_{k+1,k+1}\Phi^{-1})\nonumber\\
(\p\pbar X)_k^v&=&\gamma_{1,1}^{-1}\left(\p
\overline{Q}_{j+1}+{1\over 2}\tr(\pbar
X[v_{j+1},A_{\lambda}])-\alpha_{2,1}Q_{j+1}-\alpha_{2,2}\overline{Q}_{j+1}\right)\nonumber\\
(\p\pbar X)_k^u&=&\gamma_{1,1}^{-1}\left(\p
\overline{J}_{j+1}+{1\over 2}\tr(\pbar
X[u_{j+1},A_{\lambda}])-\alpha_{2,1}J_{j+1}-\alpha_{2,2}\overline{J}_{j+1}\right)\nonumber\\
(\p\pbar X)_k^d&=&-i\tr(\p\pbar X\Phi
E_{k+1,k+1}\Phi^{-1})\nonumber\\
\alpha_{1,1}&=&(\det G)^{-1}({1\over 2}\p g_{11}\overline{g}_{11}-\p g_{12}g_{12}+{1\over 2}\pbar g_{11}g_{12})\nonumber\\
\alpha_{1,2}&=&(\det G)^{-1}(g_{11}\p g_{12}-{1\over 2}g_{11}\pbar
g_{11}-{1\over 2}g_{12}\p g_{11})\nonumber \\
\alpha_{2,1}&=&(2\det G)^{-1}(\pbar g_{11}\overline{g}_{11}-\p \overline{g}_{11}g_{12})\nonumber\\
\alpha_{2,2}&=&(2\det G)^{-1}(g_{11}\p\overline{g}_{11}-\pbar
{g}_{11}g_{12})\nonumber\\
\det
G&=&g_{11}g_{22}-g_{12}^2=\left|\sum_{k,l=1}^{N+1}(y_k-y_l)^2(\Xi_1)_{kl}(\Xi_1)_{lk}\right|^2\nonumber\\
&-&\left(\sum_{k,l=1}^{N+1}(y_k-y_l)^2(\Xi_1)_{kl}(\Xi_2)_{lk}\right)^2\nonumber
\end{eqnarray*}
The mean curvature and the Gaussian curvature are then
\begin{eqnarray*}
H&=&(\det G)^{-1}\left(g_{11}(\pbar^2 X)^{\perp}-2g_{12}(\p\pbar
X)+g_{22}(\p^2X)^{\perp}\right)\\
K&=&(2\det G)^{-1}\Bigg(\p\left({1\over{\det G}}({g_{12}}\pbar\ln
g_{11}-\p g_{22})\right)\\
&+&\pbar \left({1\over {\det G}}(2\p g_{12}-\pbar
g_{11}-g_{12}\p\ln g_{11})\right)\Bigg)
\end{eqnarray*}
and the Willmore functional is
\begin{eqnarray*}
W=\int_{\Omega}|H|^2\sqrt{\det G}\d\xi\d\overline{\xi}
\end{eqnarray*}
\section{Explicit construction in the ${\frak su}(2)$
case}\label{se:su2}\setcounter{equation}{0}\setcounter{theorem}{0} \setcounter{definition}{0}
\setcounter{example}{0}

In general, it is difficult to obtain formula for the
matrix entries $L_m$ with the Drinfeld-Sokolov iteration process \cite{DS}.
In the $2\times 2$ case, however, one can derive general formula
for the matrix entries of $L_m$ explicitly in terms of the
solution $u_0$ by using the compatibility conditions
(\ref{eq:lax})
\begin{eqnarray*}
\p L_m&=&[L_m,A_{\lambda}]\\
\overline{\p}L_m&=&[L_m,B_{\lambda}]
\end{eqnarray*}
In this case, the matrices $A(\lambda)$ and $B(\lambda)$ take the
form
\begin{eqnarray*}
A_{\lambda}=\pmatrix{\p u_0&\lambda e^{2u_0}\cr
                     e^{-2u_0}&-\p u_0\cr}, \quad
                     B_{\lambda}=-\rho(A_{\lambda})
\end{eqnarray*}
Let the entries of $L_m$ be the following
\begin{eqnarray*}
L_m=\pmatrix{P(\lambda)&Q(\lambda)\cr
             R(\lambda)&-P(\lambda)\cr}
\end{eqnarray*}
which are polynomials in $\lambda$
\begin{eqnarray*}
P(z)&=&\sum_{i=0}^{m}P_i\lambda^i\\
Q(z)&=&\sum_{i=0}^{m}Q_i\lambda^i\\
R(z)&=&\sum_{i=0}^{m}R_i\lambda^i
\end{eqnarray*}
By comparing the coefficients of $\lambda$ in
\begin{eqnarray*}
\p L_m=[L_m,A_{\lambda}]
\end{eqnarray*}
we find that $P$, $Q$ and $R$ satisfy
\begin{eqnarray*}
\p P(\lambda)&=&Q(\lambda)e^{-2u_0}-\lambda R(\lambda)e^{2u_0}\\
\p Q(\lambda)&=&2(P(\lambda)\lambda e^{2u_0}-Q(\lambda)\p u_0)\\
\p R(z)&=&2(-P(\lambda)e^{-2u_0}+R(\lambda)\p u_0)
\end{eqnarray*}
We first write $Q=e^{-2u_0}Q^{\prime}$ and $R=e^{2u_0}R^{\prime}$.
This will simplify the last two equations
\begin{eqnarray}\label{eq:coef}
\p P(\lambda)&=&Q^{\prime}(\lambda)e^{-4u_0}-\lambda R^{\prime}(\lambda)e^{4u_0}\nonumber\\
\p Q^{\prime}(\lambda)&=&2e^{4u_0}P(\lambda)\lambda\\
\p R^{\prime}(\lambda)&=&-2e^{-4u_0}P(\lambda)\nonumber
\end{eqnarray}
By eliminating $P$ from the above equations, we can express $\p
Q^{\prime}$ in terms of $R^{\prime}$, then by eliminating
$R^{\prime}$ from the last two equations, we obtain
\begin{eqnarray*}
\lambda(16\p u_0+4\p)R^{\prime}=\left(\p^3+12\p u_0\p^2+(32(\p
u_0)^2+4\p^2u_0)\p\right)R^{\prime}
\end{eqnarray*}
We now do another change of variable
$R^{\prime}=e^{-4u_0}R^{\prime\prime}$ to obtain
\begin{eqnarray*}
4\lambda\p R^{\prime\prime}=\left(\p^3+(-8\p^2u_0-16(\p
u_0)^2)\p+(-16\p^2u_0\p u_0-4\p^3u_0)\right)R^{\prime\prime}
\end{eqnarray*}
We now denote $\Lie_n$ by the following recursion operator
\begin{eqnarray}\label{eq:Lie}
\p\Lie_n[u_0]&=&{1\over 4}\Big(\p^3+(-8\p^2u_0-16(\p
u_0)^2)\p\\&+& (-16\p^2u_0\p u_0-4\p^3u_0)\Big)\Lie_{n-1}[u_0]\nonumber\\
&=&{1\over
4}\left(\p^3-8(\p^2e^{2u_0})\p-4(\p^3e^{2u_0})\right)\Lie_{n-1}[u_0]
\nonumber\\
\Lie_0[u_0]&=&1,\quad\Lie_k[u_0]=0, \quad k<0\nonumber
\end{eqnarray}
for example, the first few terms are
\begin{eqnarray*}
\Lie_0[u_0]&=&1\\
\Lie_1[u_0]&=&-\p^2e^{2u_0}\\
\Lie_2[u_0]&=&-{1\over 16}\p^4e^{2u_0}+{9\over 8}\p^2e^{2u_0}
\end{eqnarray*}
We can now express $P$, $Q$ and $R$ in terms of $\Lie_n$
\begin{eqnarray}\label{eq:coeffin}
P&=&\sum_{i=0}^m\left(2\p u_0\Lie_{m-i-1}[u_0]-{1\over 2}\p\Lie_{m-i-1}[u_0]\right)\lambda^i\nonumber\\
Q&=&\sum_{i=0}^m\Lie_{m-i}[-u_0]e^{2u_0}\lambda^i\\
R&=&\sum_{i=0}^m\Lie_{m-i-1}[u_0]e^{-2u_0}\lambda^i\nonumber
\end{eqnarray}
We can obtain the coefficients of $\rho(L_m)$ by
\begin{eqnarray*}
\tilde{P}=\overline{P},\quad \tilde{Q}=\overline{R},\quad
\tilde{R}=\overline{Q}, \quad |\lambda|=1
\end{eqnarray*}
where $\tilde{P}$, $\tilde{Q}$, $\tilde{R}$ are the entries of the
matrix $\rho(L_m)$.

For $m=2$, the entries are
\begin{eqnarray*}
P&=&2\p u_0\lambda+{1\over 2}\p^3u_0-4(\p u_0)^3\\
Q&=&e^{2u_0}\left(\lambda^2-\p^2e^{-2u_0}\lambda-{1\over 16}\p^4e^{-2u_0}+{9\over 8}\p^2e^{-2u_0}\right)\\
R&=&e^{-2u_0}\left(\lambda-\p^2e^{2u_0}\right)
\end{eqnarray*}
We can now substitute the theta function solution (\ref{eq:sol})
for $u_0$ (as defined in proposition \ref{pro:sol}) and express
the immersion
\begin{eqnarray*}
X=i(L_m+\rho(L_m))
\end{eqnarray*}
in terms of theta function. In the case of $m=2$, this gives the
matrix entries as follows
\begin{eqnarray}\label{eq:2x2ex}
P&=&\p\log\Theta\lambda+{1\over 4}\p^3\log\Theta-{1\over 2}(\p\log\Theta)^3\nonumber\\
Q&=&\Theta\left(\lambda^2-\p^2\Theta^{-1}\lambda-{1\over 16}\p^4\Theta^{-1}+{9\over 8}\p^2\Theta^{-1}\right)\\
R&=&\Theta^{-1}\left(\lambda-\p^2\Theta\right)\nonumber
\end{eqnarray}
where $\Theta$ is the theta function expression on the right hand
side of (\ref{eq:sol}).

The components of the induced metric of the surface are the
following
\begin{eqnarray*}
g_{11}&=&-{1\over 2}\tr(\p X\p X)=-(\p
P+\p\overline{P})^2-{1\over 2}(\p(Q+\overline{R}))(\p(R+\overline{Q}))\\
g_{12}&=&-{1\over 2}\tr(\p X\pbar X)=-|(\p P+\p\overline{P})|^2-{1\over 2}(\p(Q+\tilde{Q}))(\pbar(R+\tilde{R}))\\&-&{1\over 2}(\p(R+\tilde{R}))(\pbar(Q+\tilde{Q}))\\
      &=&-|(\p P+\p\overline{P})|^2-{1\over 2}\left(|\p Q|^2+|\pbar Q|^2+|\p
      R|^2+|\pbar R|^2\right)\\&-&{1\over 2}\left({\bf Re}\p R\pbar Q+{\bf Re}\p Q\pbar R\right)\\
g_{22}&=&\overline{g}_{11}
\end{eqnarray*}
In the ${\frak su}(2)$ case, since we have
\begin{eqnarray*}
\tr(X[A,X])=\tr(X[B,X])=0
\end{eqnarray*}
we see that $\sqrt{-2\tr(X^2)^{-1}}X$ is a unit normal vector to
the surface. We have the following induced metric on the
surface ${\cal F}$
\begin{eqnarray}\label{eq:metric}
g_{11}&=&-{1\over 2}\tr(X\p^2X)=-{1\over 2}\tr(X[[X,A],A])=-{1\over 2}\tr([X,A]^2)\nonumber\\
g_{12}&=&-{1\over 2}\tr(X\p\pbar X)=-{1\over 2}\tr(X[[X,B],A])={1\over 2}\tr([B,X][X,A])\\
g_{22}&=&-{1\over 2}\tr(X\pbar^2X)=-{1\over
2}\tr(X[[X,B],B])=-{1\over 2}\tr([B,X]^2)\nonumber
\end{eqnarray}
We can now use these relations to compute the second fundamental
form and obtain the proportionality of the first and second
fundamental forms
\begin{eqnarray*}
II=\sqrt{-2\tr(X^2)^{-1}}(g_{11}\d\xi\d\xi+2g_{12}\d\xi\d\overline{\xi}+g_{22}\d\overline{\xi}\d\overline{\xi})
=(\sqrt{-2\tr(X^2)^{-1}})I
\end{eqnarray*}
Hence the mean curvature $H$ and Gauss curvature $K$ of the
surface ${\cal F}$ are constants
\begin{eqnarray*}
H&=&-(\det(G)\tr(X^2))^{-1}(2g_{11}g_{22}+2g_{12}^2)=-2\sqrt{2\over\tr(X^2)}\\
K&=&4{{\det(G)\tr(X^2)^{-2}}\over{\det(G)}}=-2\tr(X^2)^{-1}\\
\det G&=&-{1\over
2}\left(\tr([B,X][X,A])^2+\tr([B,X]^2)\tr([X,A]^2)\right)
\end{eqnarray*}
This is not surprising as $\tr(X^2)$ is a constant and hence the
surfaces are all mapped into the sphere $S^2$ of radius
$\tr(X^2)$. This means that different Toda lattice data $u_0$
correspond to different parametrizations of the same surface which
is a subset of a sphere in $\mathbb{R}^3$.

Finally we may determine formally a moving frame of the surface
${\cal F}$ and write the structural equations, namely the
Gauss-Weingarten equations and the Gauss-Codazzi equations. Let us
assume that $u_0$ is a solution of the Toda lattice equations
(\ref{eq:toda}) such that $\det G\neq 0$ in some open neighborhood
of a regular point $(\xi,\overline{\xi})$ in
$\Omega\in\mathbb{C}$. Suppose also that the surface ${\cal F}$
given by (\ref{eq:immersion}) is associated with these structural
equations. The surface ${\cal F}$ can be described by the moving
frame
\begin{eqnarray*}
\eta =(\eta_1=\p X, \eta_2=\pbar X, \eta_3,\ldots,\eta_{N^2-1})^T
\end{eqnarray*}
where $\eta_3,\ldots,\eta_{N^2-1}$ are real normal vectors to the
surface ${\cal F}$. The components of the moving frame $\eta$ can
be written in terms of $N\times N$ skew-hermitian matrices
satisfying the following normalization conditions
\begin{eqnarray*}
(\p X,\p X)&=&g_{11},\quad (\p X,\pbar X)=g_{12},\quad (\pbar
X,\pbar X)=g_{22}\\
(\p X,\eta_k)&=&(\pbar X,\eta_k)=0,\quad
(\eta_j,\eta_k)=(\eta_j,\eta_j)\delta_{jk},\quad
j,k=3,\ldots,N^2-1
\end{eqnarray*}
In the ${\frak su}(2)$ case, such a moving frame on ${\cal F}$
consists of the the vectors $\p X$, $\pbar X$ and $X$.

Let the frame $\eta=(\p X, \pbar X, X)^T$, we have the following
Gauss-Weingarten equations
\begin{eqnarray*}
\p \eta&=&M_1\eta,\quad \pbar \eta=M_2\eta\\
M_1&=&\pmatrix{\alpha_{1,1}&\alpha_{1,2}&2g_{11}\cr
\alpha_{2,1}&\alpha_{2,2}&2g_{12}\cr 1&0&0\cr} \\
M_2&=&\pmatrix{\alpha_{2,1}&\alpha_{2,2}&2g_{12}\cr
\overline{\alpha}_{1,2}&\overline{\alpha}_{1,1}&2\overline{g}_{11}\cr
0&1&0\cr}
\end{eqnarray*}
where $\alpha_{i,j}$ are given by the following
\begin{eqnarray*}
\alpha_{1,1}&=&(\det G)^{-1}({1\over 2}\p g_{11}\overline{g}_{11}-\p g_{12}g_{12}+{1\over 2}\pbar g_{11}g_{12})\nonumber\\
\alpha_{1,2}&=&(\det G)^{-1}(g_{11}\p g_{12}-{1\over 2}g_{11}\pbar
g_{11}-{1\over 2}g_{12}\p g_{11})\nonumber \\
\alpha_{2,1}&=&(2\det G)^{-1}(\pbar g_{11}\overline{g}_{11}-\p \overline{g}_{11}g_{12})\nonumber\\
\alpha_{2,2}&=&(2\det G)^{-1}(g_{11}\p\overline{g}_{11}-\pbar
{g}_{11}g_{12})\nonumber
\end{eqnarray*}
which are the usual Christoffel symbols of second kind.

The coefficients $\alpha_{1,1}$ and $\alpha_{1,2}$ can be computed
by solving the following linear equations
\begin{eqnarray}\label{eq:lineareq}
\tr\left((\p^2X-\alpha_{1,1}\p X-\alpha_{1,2}\pbar X)\p X\right)&=&0 \nonumber\\
\tr\left((\p^2X-\alpha_{1,1}\p X-\alpha_{1,2}\pbar X)\pbar
X\right)&=&0
\end{eqnarray}
Since from (\ref{eq:lineareq}), we have
\begin{eqnarray*}
\p^2X=\alpha_{1,1}\p X+\alpha_{1,2}\pbar X+\gamma X
\end{eqnarray*}
for some $\gamma$, a solution $\alpha_{1,1}$ and $\alpha_{1,2}$ to the (\ref{eq:lineareq})
would give the $\p X$ and $\pbar X$ components of
$\p^2X$. To determine $\gamma$, note that since $X$ is
orthogonal to the plane spanned by $\p X$ and $\pbar X$, the
coefficient $\gamma$ is just the inner product between $\p^2X$ and
$X$
\begin{eqnarray*}
\gamma={1\over 2}\tr(\p^2X X)=g_{11}
\end{eqnarray*}
which follows from (\ref{eq:metric}). The other coefficients can
be computed similarly.

\section{Conclusions}\setcounter{equation}{0}\setcounter{theorem}{0} \setcounter{definition}{0}
\setcounter{example}{0}
In this paper we have studied the ${\frak su}(N)$ algebraic
description of surfaces obtained from theta function solutions of
the 2D Toda lattice. We have derived the generalized Weierstrass
formula for immersion which is expressed in terms of the Lax
operators. It allows us to construct surfaces immersed in ${\frak
su}(N)$ algebra and show that all these surfaces are mapped into hyperspheres
in $\mathbb{R}^{N^2-1}$. It has been proved to be effective, since we were
able to recover easily the known results associated to the ${\frak
su}(2)$.

The structural equations of the 2-dimensional surfaces, their
moving frames, the first and second fundamental forms, the
Gaussian curvature, the mean curvature vector have been expressed
in terms of the Baker function of the Toda lattice. The
implementation of these theoretical results have been applied to
the ${\frak su}(2)$ case, leading to the negative constant
Gaussian curvature, which means that these surfaces are subsets of
a sphere in $\mathbb{R}^3$.

The approach presented here is limited to theta function solutions
of the 2D Toda lattice. The question arises as to whether this
approach can be extended to integrable systems (for example to the
$\mathbb{CP}^N$ models in (2+1)-dimensions or non-Abelian gauge
theories) and if these models can lead to the construction of
other classes of submanifolds immersed in multi-dimensional
spaces. Other aspects of the method worth investigating are the
completeness of solutions, the compactness of surfaces and the
stability of the resulting surfaces.

The geometrical analysis of surfaces and their deformations under
various types of dynamics have generated a great deal of interest
and activities in several mathematical and physical areas of
research as well. In particular, it is worth noting that we know
in certain cases the analytic description of surfaces in a
physical system for which analytic models are not yet fully
developed. However, by using our approach we can in some cases
select an appropriate Toda lattice model corresponding to the
given Weierstrass representation and characterize the classes of
equations describing the physical phenomena under investigation.
This approach was attempted successfully for the Weierstrass
representation associated to sine-Gordon model, for surfaces with
negative Gaussian curvature in three-dimensional Euclidean space
\cite{BK}, \cite{RS}, but not to our knowledge for
multi-dimensional spaces.

\end{document}